\documentclass[11pt,reqno]{amsart}

\usepackage{amsmath,amssymb, mathabx}
\usepackage{graphicx}
\usepackage[pagebackref, colorlinks = true, linkcolor = blue, urlcolor  = blue, citecolor = red]{hyperref}
\usepackage{xcolor}

\usepackage[margin=1in]{geometry}

\usepackage{caption} 
\usepackage{subcaption}

\newcommand{\eq}[1]{\begin{align}\begin{aligned}#1\end{aligned}\end{align}}

\newtheorem{theorem}{Theorem}[section]
\newtheorem*{theorem*}{Theorem}
\newtheorem{definition}[theorem]{Definition}
\newtheorem{proposition}[theorem]{Proposition}
\newtheorem{corollary}[theorem]{Corollary}
\newtheorem{lemma}[theorem]{Lemma}
\newtheorem{remark}[theorem]{Remark}

\DeclareMathOperator{\trace}{Tr} 
\DeclareMathOperator{\id}{id} 

\DeclareMathOperator{\inter}{Int}
\DeclareMathOperator{\nc}{NC}
\newcommand{\kr}{\mathrm{Kr}}
\newcommand{\Boole}{\mathrm{boole}}
\newcommand{\free}{\mathrm{free}}

\usepackage{layout}
\date{\today}

\title{Generating series and matrix models for \\ meandric systems with one shallow side}

\author{Motohisa Fukuda}
\address{MF: Yamagata University, 1-4-12 Kojirakawa, Yamagata, 990-8560 Japan}
\email{fukuda@sci.kj.yamagata-u.ac.jp}

\author{Ion Nechita}
\address{IN: Laboratoire de Physique Th\'eorique, Universit\'e de Toulouse, CNRS, UPS, France}
\email{nechita@irsamc.ups-tlse.fr}

\subjclass[2000]{}

\begin{document}
\maketitle
\markright{\MakeUppercase{Generating series and matrix models for meandric systems with one shallow side}}

\begin{abstract}
In this article, we investigate meandric systems having one shallow side: the arch configuration on that side has depth at most two. This class of meandric systems was introduced and extensively examined by I.~P.~Goulden, A.~ Nica, and D.~Puder in \cite{goulden2020asymptotics}. Shallow arch configurations are in bijection with the set of interval partitions. We study meandric systems by using moment-cumulant transforms for non-crossing and interval partitions, corresponding to the notions of free and boolean independence, respectively, in non-commutative probability. We obtain formulas for the generating series of different classes of meandric systems with one shallow side, by explicitly enumerating the simpler, irreducible objects. In addition, we propose random matrix models for the corresponding meandric polynomials, which can be described in the language of quantum information theory, in particular that of quantum channels. 
\end{abstract}

\tableofcontents

\section{Introduction}

Meanders are fundamental combinatorial objects of great complexity, defined by a simple 
non-crossing closed curve intersecting a reference line at $2n$ points. Their enumeration (as a function of $n$) is an important open problem in combinatorics \cite{albert2005bounds}. There is a large theoretical body of work dealing with the combinatorics of meanders, see \cite{lando1993plane,difrancesco1997meander}.

Mathematically, meandric systems are generalizations of meanders consisting of two \emph{arch configurations}, one on top and the other one on the bottom of the reference line. An arch configuration corresponds precisely to a non-crossing pairing of the coordinate set $\{1,2, \ldots, 2n\}$. It is this connection to the theory of non-crossing partitions that had been put forward by A.~Nica in \cite{nica2016free}, starting the study of meandric systems with the help of tools from free probability theory \cite{voiculescu1985symmetries,mingo2017free}. This line of work has been pursued further, with new results about semi-meanders \cite{nica2018operator}, or meandric systems with large number of loops \cite{fukuda2019enumerating}.

An important result was obtained by I.P.~Goulden, A.~Nica, and D.~Puder in \cite{goulden2020asymptotics}, where a particular subclass of meandric systems was described combinatorially: the authors studied meandric systems where the arch configurations on top of the reference line correspond to interval partitions; the authors named such meandric systems
\emph{shallow top meanders}. Due to the simpler combinatorial structure of the top arch configurations, shallow top meanders are tractable 
enough to provide interesting lower bounds on the total number of meanders. 

Our work drew most of its inspiration from \cite{goulden2020asymptotics}, but tackles the enumeration of special classes of meandric systems in a systematic way, employing tools from non-commutative probability theory. Our main insight is to reduce the enumeration of meandric systems to that of a simpler class of objects, sometimes called ``irreducible'' (see \cite{beissinger1985enumeration} for the general flavor in combinatorics). If the initial class of meandric systems corresponds to the \emph{moments} of some non-commutative distribution, the simpler meandric systems correspond to its \emph{cumulants}. The type of cumulants involved depends on the structure of the initial meandric systems: general non-crossing partitions yield free cumulants, while interval partitions boolean cumulants. Once the probabilistic machinery is applied, we can then directly enumerate the simpler combinatorial objects
and in theory the initial, allowing us to treat several situations in a unified manner. 

Historically, meandric systems were also studied using methods from random matrix theory. P.~Di Francesco and his collaborators developed several such models in \cite{difrancesco1997meander,di2001matrix}. Later, an intriguing connection to the theory of quantum information theory was put forward in \cite{fukuda2013partial}. We provide at the end of this paper several new matrix models for the various classes of meandric systems we consider, which also fall in the field of quantum information. Indeed, we show that the meandric polynomial is equal to the asymptotic moments of the output state of a tensor product of completely positive maps, acting on the maximally entangled state. The choice of completely positive maps depends on the type of partitions on the bottom side one considers: random Gaussian channels for general non-crossing partitions and a depolarizing channel for interval partitions. These models are conceptually simpler than the past ones, and allow us to treat the different subsets of meanders in a unified manner. 

Our paper is organized as follows. Section \ref{sec:combinatorics} contains the main definitions and tools from the combinatorial theory of permutations and meanders. In Section \ref{sec:transforms} we recall the basic tools from boolean and free probability theory used in this work. The following three sections contain the main body of the paper, dealing with three different classes of meandric systems: thin (both shallow top and shallow bottom) meandric systems in Section \ref{sec:thin}, shallow top meandric systems in Section \ref{sec:shallow-top}, and shallow-top semi-meanders in Section \ref{sec:shallow-top-semi}. Finally, random matrix models are discussed in Section \ref{sec:RMT}.

\section{Combinatorial aspects of meandric systems} \label{sec:combinatorics}

\subsection{Basics of non-crossing partitions and permutations}\label{sec:NC-permutations}

This section contains the necessary definitions and properties of the combinatorial objects meandric systems are built on, which are mainly non-crossing and interval partitions. We refer the reader to \cite{biane1997some} or \cite{nica2006lectures} for more details.

We denote by $\mathcal S_n$ the group of permutations of $n$ symbols. For a permutation $\alpha \in S_n$, we denote by $\|\alpha\|$ its length: $\|\alpha\|$ is the minimal number $m$ of transpositions $\tau_1, \ldots ,\tau_m$ which multiply to $\alpha$:
$$\|\alpha\| := \min\{m \geq 0 \, : \, \exists \, \tau_1, \ldots, \tau_m \in \mathcal S_n \text{ transpositions s.t.~}
 \alpha = \tau_1 \cdots \tau_m\}.$$
The length $\| \cdot \|$ endows the symmetric group $\mathcal S_n$ with a metric structure, by defining $d(\alpha, \beta) = \|\alpha^{-1}\beta\|$.
The following relation between the number of cycles $\#\alpha$ of a permutation and its length is crucial to us:
\eq{
\|\alpha\| + \#(\alpha) = n.
}
Both statistics $\#(\cdot)$ and $\| \cdot \|$ are constant on conjugation classes, hence the following relations hold:
$$\|\alpha\| = \|\alpha^{-1}\|\qquad \text{ and } \qquad \|\alpha \beta \| = \|\beta \alpha\|.$$ 

\begin{figure}[htbp!]
\centering
\begin{minipage}{.49\textwidth}
  \centering
  \includegraphics[scale=.6]{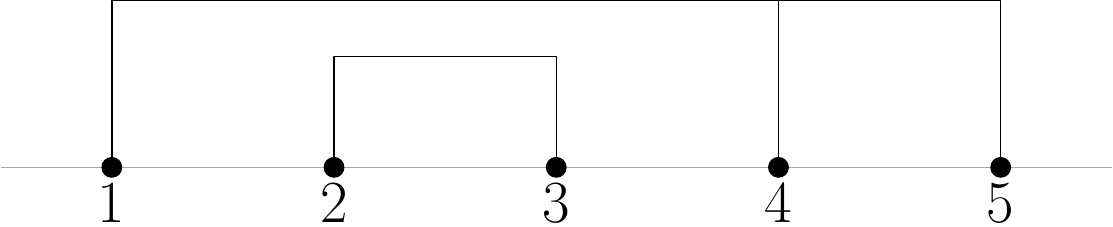}
  \captionof{figure}{A non-crossing partition.}
  \label{fig:NC1}
\end{minipage}\begin{minipage}{.49\textwidth}
  \centering
  \includegraphics[scale=.6]{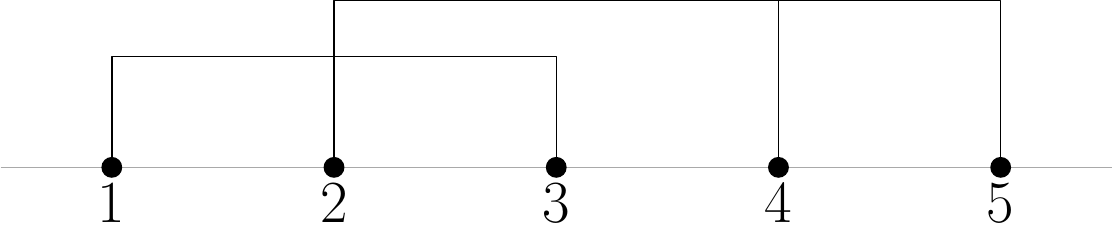}
  \captionof{figure}{A crossing partition.}
  \label{fig:NC2}
\end{minipage}
\end{figure}

Let us now introduce the different classes of partitions which will be of interest to us. A partition $B_1 \sqcup B_2 \sqcup \cdots \sqcup B_m = \{1,2,\ldots,n\} =:[n]$ is called \emph{non-crossing} if its blocks $B_k$ do not cross:
there do not exist distinct blocks $B_i, B_j$ and $a,b \in B_i$ and $c,d \in B_j$ such that $a<c<b<d$. The partition $\{1,4,5\} \sqcup \{2,3\}$ of $[5]$ is non-crossing, see Figure \ref{fig:NC1}. The partition $\{ 1,3\}\sqcup\{2,4,5\}$ on the other hand is crossing, see Figure \ref{fig:NC2}. The set of non-crossing partitions of $[n]$ is denoted by $\nc(n)$ or $\nc(1,2,\ldots,n)$, if we want to emphasize the underlying set. The subset of non-crossing partitions consisting of \emph{pairings} (i.e.~all the blocks have size two) is denoted by $\nc_2(n)$; in this case, $n$ must obviously be even. Finally, the subset of \emph{interval} partitions, denoted by $\inter(n)$ consists of (non-crossing) partitions having blocks made of consecutive integers. We have 
$$|\nc(n)| = \mathrm{Cat}_n = \frac{1}{n+1}\binom{2n}{n} \qquad \text{ and } \qquad |\inter(n)| = 2^{n-1}.$$

In many cases, it is important to identify non-crossing partitions with a class of permutations, called \emph{geodesic permutations}. This correspondence, initially observed in \cite{biane1997some} (see also \cite[Lecture 23]{nica2006lectures}) is key in many areas, and used extensively in random matrix theory for example. The bijection is defined as follows: one associates to each block of a non-crossing partition a cycle in a permutation where the elements are ordered increasingly. For example, the non-crossing partition from Figure \ref{fig:NC1} is identified to the permutation $(1,4,5)(2,3) \in \mathcal S_5$. Note that the cycles of the permutation are precisely the blocks of the non-crossing partition, with the choice of (cyclically) ordering the elements increasingly.  

Importantly, it was shown by Biane \cite{biane1997some} that geodesic permutations are characterized using the metric induced by the length function on the symmetric group: $\alpha \in \mathcal S_n$ is a geodesic permutation
if and only if it saturates the triangle inequality
$$\|\alpha\| + \|\alpha^{-1} \gamma \| = \|\gamma\| = n-1,$$
where $\gamma = (1,2, \ldots n)$ is the full-cycle permutation. Note that in this case 
$\alpha$ lies on the geodesic between the identity permutation $\id=(1)(2)\cdots(n)$ and $\gamma=(1,2,\ldots,n)$.

The set $NC(n)$ is endowed with a partial order called \emph{reversed refinement}:  $\alpha \leq \beta$ if every block of $\alpha$ is contained in a block of $\beta$. Note that this order relation is not total: for example, partitions $\{1\} \sqcup \{2,3\}$ and $\{2\} \sqcup \{1,3\}$ are not comparable. This partial order can be nicely characterized in terms of the associated geodesic permutations: $\alpha \leq \beta$ is equivalent to $\alpha$  lying on the geodesic between $\id$ and $\beta$: 
$$\|\alpha\| + \|\alpha^{-1} \beta\| = \|\beta\|.$$

\begin{figure}[htbp]
  \centering
    \includegraphics[width=0.6\textwidth]{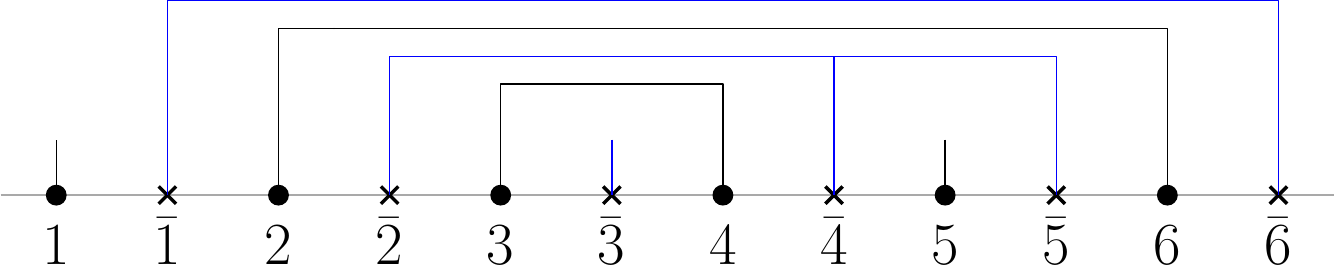}
\caption{The Kreweras complement of $\alpha = \{1\} \sqcup \{2,6\} \sqcup \{3,4\} \sqcup \{5\}$ is \textcolor{blue}{$\alpha^\kr = \{1,6\} \sqcup \{2,4,5\} \sqcup \{3\}$}.}
 \label{fig:Krew2}
\end{figure}

Let us now discuss the important notion of \emph{Kreweras complement} for non-crossing partitions. The Kreweras complement is an order reversing involution $\alpha \mapsto \alpha^\kr$ of $\nc(n)$, defined in the following way \cite[Definition 9.21]{nica2006lectures}. 
First, double the elements of the basis set to obtain $\{1, \bar 1,  2,\bar2 \ldots,  n,\bar  n\}$
and then  consider $\alpha^\kr \in \nc(\bar 1,\bar 2, \ldots, \bar n) \cong \nc(n)$ be 
the largest non-crossing partition such that $\alpha \sqcup \alpha^\kr$ is still a non-crossing partition on $\{1,\bar 1 ,2, \bar2 \ldots ,  n, \bar n\}$. This operation is best explained by an example, see Figure \ref{fig:Krew2}: for $\alpha = (2,6)(3,4)$ we have $\alpha^\kr=(1,6)(2,4,5)$. The extremal elements in $\nc(n)$ are swapped: $\id^\kr=\gamma$ and $\gamma^\kr=\id$. In the language of geodesic permutations, given a geodesic permutation $\id-\alpha-\gamma$, 
the Kreweras complement of $\alpha$ corresponds to the permutation $\alpha^\kr \in \mathcal S_n$ defined as
\eq{\label{eq:kreweras_rewritten}
\alpha^\mathrm{Kr} = \alpha^{-1} \pi
}
see \cite[Remark 23.24]{nica2006lectures} for details. Importantly, for $\alpha \in \nc(n)$, we have
\eq{\label{eq:kreweras_norm}
\| \alpha \| + \|\alpha^\kr \| = n-1.
}

Finally, let us discuss the bijection between $\nc(n)$ and $\nc_2(2n)$, called \emph{fattening} (note that both sets are counted by the Catalan numbers). For a given non-crossing partition $\alpha \in \nc(n)$, we consider two points $i_{-}$ and $i_{+}$ for both sides of each $i \in \{1,\ldots,n\}$,
left and right respectively, doubling in this way the index set. We associate to $\alpha$ the following pairing: connect $i_{+}$ and $j_{-}$ if $\alpha(i) =j$, where $\alpha$ is seen now as a permutation. It can be shown that the pair partition obtained in this way is non-crossing, see \cite[Lecture 9]{nica2006lectures} for the details.

\subsection{Loops in meandric systems}
\begin{figure}[htbp!] 
  \centering
    \includegraphics[width=1\textwidth]{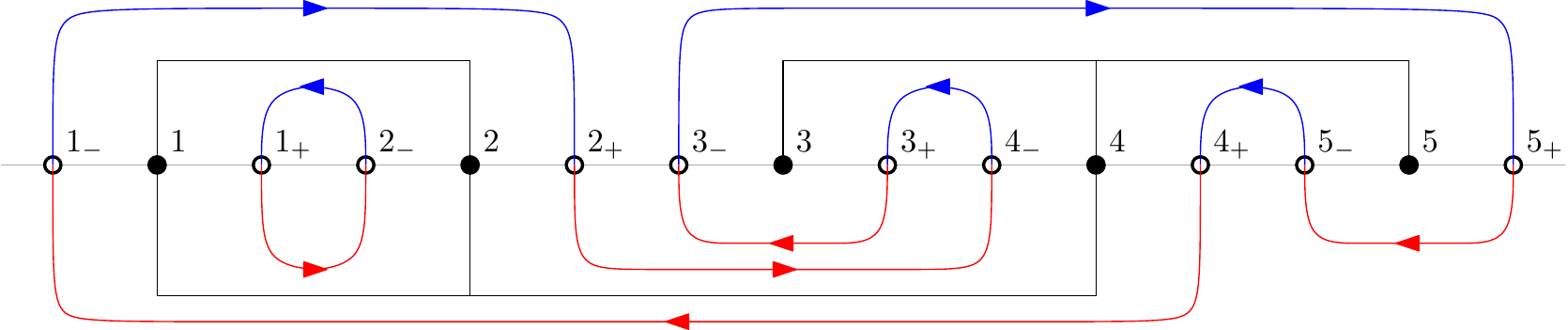}
\caption{A meandric system generated by the geodesic permutations (in black) $\alpha = (1,2)(3,4,5)$ and $\beta = (1,2,4)(3)(5)$.
The two loops (blue and red together) with arrows are formed by the action of $\alpha^{-1}\beta$.}
 \label{fig:MtoP}
\end{figure}

As discussed in the introduction, there has been a lot of interest in counting meandric systems with respect to their number of connected components, which we call \emph{loops}. In this paper, we shall regard meandric systems as pairs of non-crossing partitions (or geodesic permutations). This point of view is best explained with an example, see Figure  \ref{fig:MtoP}. In this figure, the meandric system is made of the blue and red arches, connecting the points $\{i_{\pm}\}_{i \in [5]}$. The blue (resp.~red) 
arches on top (resp.~bottom) on the reference line are associated to non-crossing pairings (called arch configurations in \cite{difrancesco1997meander}), which, in turn, are in bijection to the non-crossing partitions connecting the points $\{i\}_{i \in [5]}$ displayed in black. In the figure, the black line above and below the reference line correspond to non-crossing partitions $\alpha = (1,2)(3,4,5)$ and $\beta = (1,2,4)(3)(5)$,
respectively. The blue and red lines are fattenings of those permutations, which are non-crossing pairings generating the meandric system. In this example, the number of loops in this meandric system is $2$. Remarkably, it can be calculated by
\eq{
\#(\alpha^{-1} \beta) = \#( (1,2)(5,4,3)\circ(1,2,4)(3)(5)) = \#((1) (2,3,5,4)) = 2
}
To see this, one can follow the arrows in the figure to count the number of loops. 
In addition, note that in this example the top side is shallow while the bottom not. 

In short, graphically, two permutations over and under the straight lines give structural lines. ``Fattening'' them (or drawing new lines both side of those lines) gives loops of the meandric system.  
We state this property in general in the following proposition. One can refer to \cite[Section 3]{nica2016free} or \cite[Proposition 3.1]{fukuda2019enumerating} for the proof.
 
\begin{proposition}\label{proposition:loops}
Suppose a meandric system on $2n$ points is generated by $\alpha, \beta \in \nc(n)$. Then, the number of loops of the meandric system is $\#(\alpha^{-1}\beta)$, the number of cycles of the permutation $\alpha^{-1}\beta$.
\end{proposition}

The result above is crucial to our work, since it allows us to relate the problem of counting loops of meandric systems to a combinatorial problem on (special subsets of) the symmetric group. We shall also need the following lemma, showing that the Kreweras complement operation does not change the statistics of the systems. 
\begin{proposition}\label{proposition:kreweras_loop}
For $\alpha, \beta \in \nc(n)$ we have 
\eq{\label{eq:kreweras-loops}
\# ((\alpha^{\kr})^{-1}\beta^{\kr}) = \# (\alpha^{-1}\beta)
}
\end{proposition}
\begin{proof}
By using the property of Kreweras complement we have
\eq{
(\alpha^{\kr})^{-1}\beta^{\kr} = (\alpha^{-1} \gamma_n)^{-1} (\beta^{-1} \gamma_n)
= \gamma^{-1} \alpha \beta^{-1} \gamma_n \ .
}
Since $\#(\cdot)$ is a class function, we have proved the claim. 
\end{proof}

\section{Free and boolean transformations}\label{sec:transforms}

Our approach for the enumeration of special subsets of meandric systems is based on the theory of non-commutative probability theory. More precisely, using the various notions of independence existing in the non-commutative setting, we decompose meandric systems in irreducible components via the corresponding moment-cumulant formulas, which we then proceed to enumerate. In this section we gather the relevant facts and formulas from the theory of free and boolean independence, as well as some related technical combinatorial results that will be used in the later sections. 
\subsection{Basics}
Below we discuss structures of the non-crossing partitions and the interval partitions, and associated transforms. Although these notions stem from various notions of non-commutative independence, we shall not make use of the probabilistic interpretations, and focus on the combinatorics. Concretely, we shall use these transforms to relate the generating series of some combinatorial class (encoded by the moments of some non-commutative distribution) to the generating series of a simpler class (encoded by the cumulants of some type). In the free probability theory, one can define inevitable transforms called
\emph{moment-cumulant formula}: for a lattice  $L(n) \in\{ \inter(n) ,\nc(n)\}$ it holds that 
\eq{
\varphi(a_1\cdots a_n) = \sum_{\sigma\in L(n)} \kappa_{\sigma}[a_1, \ldots,a_n] 
}
where $\varphi(\cdot)$ and $\kappa(\cdot)$ are respectively moment and cumulant functionals (depending on $L$), and $a_i$'s are non-commutative random variables. 
Interested readers can refer to \cite{speicher1993boolean,lehner2004cumulants,nica2006lectures}. Now, restricting ourselves to the case $a = a_1 = \cdots = a_n$ and using the multiplicativity of  $\kappa(\cdot)$, we define the following transformations:
\begin{definition}\label{definition:moment-cumulant2}
Between sequences of numbers, 
\emph{boolean transform} $\mathcal{F}_{\Boole}$ and \emph{free transform} $\mathcal{F}_{\free}$:
\eq{
\mathcal{F}_{\cdot} : \{\kappa_n\}_{n=1}^\infty \mapsto \{m_n\}_{n=1}^\infty
}
are defined by
\eq{
m_n = \sum_{\sigma \in L(n)} \prod_{c \in \sigma} \kappa_{|c|} \ ,
}
Here, $L(n) = \inter(n)$ in the case of $\mathcal{F}_{\Boole}$ and $L(n) =  \nc(n)$ in the case of $\mathcal{F}_{\free}$, while $c \in \sigma$ are the blocks of $\sigma$.
This can be extended naturally to maps between polynomials (moment and cumulant generating functions):
\eq{
\mathcal{F}_{\cdot} : K(X) \mapsto M(X)
}
where
\eq{
M(X) = \sum_{n = 1}^\infty m_n X^n \qquad \text{ and } \qquad
K(X) = \sum_{n = 1}^\infty \kappa_n X^n \ .
}
\end{definition}

We quote a well-known property of the boolean Transform:
\begin{proposition}[Functional relation for boolean transform {\cite[Proposition 2.1]{speicher1993boolean}}]\label{proposition:Boole}
Suppose the moment and cumulant generating functions $M = M(X)$ and $K = K(X)$ are related 
through the boolean transform as in Definition \ref{definition:moment-cumulant2}: $\mathcal F_{\Boole}:K \mapsto M$. Then,
\eq{
K = \frac{M}{1+M}, \quad\text{ or }\quad M = \frac{K}{1-K} \ .
}
\end{proposition}

Next, we state a simple generalization of the moment-cumulant formula for free independence \cite[Lecture 11]{nica2006lectures} which treats the last block (i.e.~the block containing $n$ for a partition $\beta \in \nc(n)$) separately. Recall first that two generating series $K,M$ related by the free transform $\mathcal F_{\free}:K \mapsto M$ are related by the \emph{implicit} equation
$$M(X) = K(X(1+M(X)).$$

\begin{lemma}\label{lemma:free2}
For two sequences $\{h_n\}_{n=1}^\infty$ and $\{g_n\}_{n=1}^\infty$,
we have
\eq{
\sum_{n=1}^\infty X^n \sum_{ \beta \in \nc (n) } 
h_{|\beta(n)|} \prod_{c \in \beta^\prime} g_{|c|} 
=  \sum_{s=1}^\infty    h_{s} X^s
\left( 1+ \sum_{i=1}^\infty \hat g_i X^i \right)^s \ .
}
Here, $\mathcal F_{\free} : \{g_n\}_{n=1}^\infty \mapsto \{\hat g_n\}_{n=1}^\infty$
in Definition \ref{definition:moment-cumulant2}, and we used the following decomposition:
\eq{\label{eq:decomposition}
\beta = \beta^\prime \oplus \beta(n)
}
where $\beta(n)$ is the block of $\beta$ containing $n$.
\end{lemma}
\begin{proof}
Our proof is a standard computation:
\eq{
 \sum_{n=1}^\infty X^n \sum_{ \beta \in \nc (n) } 
h_{|\beta(n)|} \prod_{c \in \beta^\prime} g_{|c|} 
&= \sum_{n=1}^\infty   \sum_{s=1}^n h_{s} X^s
  \sum_{\substack{i_1 + \ldots + i_s = n-s \\ \text{with } i_j \geq 0}}  \, \prod_{j=1}^s \, \underbrace{\sum_{\beta^\prime_j \in \nc (i_j)} \, \prod_{c \in \beta^\prime_j }g_{|c|} X^{i_j}}_{(\clubsuit)}\\
&= \sum_{s=1}^\infty    h_{s} X^s
  \sum_{m=0}^\infty \, \sum_{\substack{i_1 + \ldots + i_s =m \\ \text{with } i_j \geq 0}}  \, \prod_{j=1}^s  \hat g_{i_j} X^{i_j} \\
&= \sum_{s=1}^\infty    h_{s} X^s
\left( 1+ \sum_{i=1}^\infty \hat g_i X^i \right)^s \ .
}
Here, we have $(\clubsuit) = 1$ when $i_j = 0$,
which corresponds to the Catalan number $\mathrm{Cat}_0 =1$.
\end{proof}

\subsection{Join and meet} 
On the lattice of $\nc(n)$ two important operations are defined.
One is so-called \emph{join} 
the smallest element $\gamma \in NC(n)$ such that $\gamma \geq \alpha,\beta$.
The other is so-called \emph{meet} 
the largest element $\gamma \in NC(n)$ where 
$\gamma \leq \alpha,\beta$: for $\alpha,\beta \in \nc(n)$ 
\eq{\label{definition:vee-nc}
&\text{join:}  &\alpha \vee \beta &= \min \{\gamma \in \nc(n): \alpha, \beta \leq \gamma \} \ , \\
&\text{meet:} &\alpha \wedge \beta &= \max \{\gamma \in \nc(n): \alpha, \beta \geq \gamma \}  \ .
}
Here, the smallest and the largest elements in $NC(n)$ are denoted by $0_n$ and $1_n$, such that
$0_n = (1)\cdots(n)$ and $1_n = (1,\ldots,n)$.
Note that $(\alpha \wedge \beta)^\mathrm{Kr}=\alpha^\mathrm{Kr} \vee \beta^\mathrm{Kr}$ and 
$(\alpha \vee \beta)^\mathrm{Kr}=\alpha^\mathrm{Kr} \wedge \beta^\mathrm{Kr}$.

First, we restrict the operations meet and join to $\inter(n) \subseteq \nc(n)$. To this end we denote the complement of $\inter(n)$ by
\eq{
\kr\inter(n) = \{\alpha^{\kr}: \alpha \in \inter(n) \} \ ,
}
and give:
\begin{definition}\label{definition:vee-int}
We define join and meet in $\inter(n)$: for $\alpha,\beta \in \nc(n)$ 
\eq{
\alpha \vee_{\inter} \beta &= \min \{\gamma \in \inter(n): \alpha, \beta \leq \gamma \} \ , \\
\alpha \wedge_{\kr \inter} \beta &= \max \{\gamma \in \kr \inter: \alpha, \beta \geq \gamma \}  \ .
}
\end{definition}
Similarly as before we have
\eq{\label{eq:meet-join-kr}
(\alpha \vee_{\inter} \beta)^\mathrm{Kr}=\alpha^\mathrm{Kr} \wedge_{\kr \inter}  \beta^\mathrm{Kr}
}
\begin{remark}
The notion $\alpha \vee_{\inter} \beta$ in Definition \ref{definition:vee-int} coincides with the definition of  ``interval closure'' in \cite[Definition 2.3---(11)]{arizmendi2015relations}.
\end{remark}

\begin{lemma}[Key decomposition]\label{lemma:decomposition}
We have the following identification:
\eq{
\inter(n) \times \nc(n) = 
\left\{(\sigma, \alpha, \beta) \in \inter(n) \times \inter(n) \times \nc(n):
\alpha \vee_{\inter} \beta = \sigma\right\}}
Moreover, we have the following bijective map:
for fixed $\sigma = c_1 \cdots c_m \in \inter(n)$, where $c_i$'s are blocks of $\sigma$, 
\eq{
&\left\{(\alpha, \beta) \in \inter(n) \times \nc(n): \alpha \vee_{\inter} \beta = \sigma \right\}\\
&\to \bigtimes_{i=1}^m \left\{(\alpha_i,\beta_i )\in \inter(|c_i|) \times \nc(|c_i|) : \alpha \vee_{\inter}\beta = 1_{|c_i|}\right\} \ .
}
Also, a similar one-to-one relation holds true after replacing $\nc(n)$ by $\inter(n)$.
\end{lemma}
\begin{proof}
Since the first identification is just a matter of classification,
we prove the second. First, we define the map. The condition $\sigma = \alpha \vee_{\inter} \beta$ has two implications. One is that 
we can write $\alpha = \oplus_{i=1}^m \alpha|_{c_i}$ and $\beta = \oplus_{i=1}^m \beta|_{c_i}$ where $\alpha|_{c_i} \in \inter(|c_i|)$ and $\beta|_{c_i} \in \nc(|c_i|)$, i.e. each block of $\alpha$ and $\beta$ belongs to one of $c_i$'s.
This is because $\alpha \vee_{\inter} \beta$ would be coarser otherwise. 
The other is that $\alpha|_{c_i} \vee_{\inter} \beta|_{c_i} = 1_{|c_i|}$
because $\alpha \vee_{\inter} \beta$ would be finer otherwise. 
Next, it is clear that the map is injective, because if two permutations are identical on each sub-interval, they are necessarily the same. 
Finally, to show surjectivity, take
$\alpha|_{c_i} \in \inter(|c_i|)$ and $\beta|_{c_i} \in \nc(|c_i|)$ with $\alpha|_{c_i} \vee_{\inter} \beta|_{c_i} = 1_{|c_i|}$,
and form $\alpha = \oplus_{i=1}^m \alpha|_{c_i} \in \inter(n)$ and $\beta = \oplus_{i=1}^m \beta|_{c_i} \in \nc(n)$.
The construction implies that $\alpha \vee_{\inter} \beta \leq \sigma$,
and the condition $\alpha|_{c_i} \vee_{\inter} \beta|_{c_i} = 1_{|c_i|}$ implies that $\alpha \vee_{\inter} \beta \geq \sigma$. This completes the proof.
\end{proof}

\begin{definition}\label{definition:MK}
For $L(n) = \nc(n) \text{ or } \inter(n)$, define the following sets:
\eq{
M_{n,r,a,b} &=\left\{(\alpha,\beta) \in \inter(n) \times L(n) :\|\alpha^{-1} \beta\| = r,\, \|\alpha\| =a,\,\|\beta\| =b   \right\} \\
K_{n,r,a,b}&=  \left\{(\alpha,\beta) \in \kr \inter(n) \times \kr L(n) : \right. \\
&\hspace{5mm} \left. \|\alpha^{-1} \beta\| = r,\, \|\alpha^{ -1} 1_n\| =a,\,\|\beta^{-1} 1_n\| =b,\,  \alpha\wedge_{\kr \inter} \beta = 0_n \right\}
}
and functions:
\begin{align}
\label{eq:def-M} M(X,Y,A,B) &= \sum_{n=1}^\infty m_n X^n 
\quad &\text{where} \qquad m_n &= 
\sum_{\substack{\alpha \in \inter(n) \\ \beta \in L (n)}} 
Y^{\|\alpha^{-1} \beta\|}A^{\|\alpha\|} B^{\|\beta\|} \\
\label{eq:def-K} K(X,Y,A,B) &=\sum_{n=1}^\infty \kappa_n X^n 
\quad &\text{where} \qquad \kappa_n &= 
\sum_{\substack{\alpha \in \kr \inter(n) \\ \beta \in \kr L(n) \\ \alpha \wedge_{\kr \inter} \beta= 0_n}} 
Y^{\|\alpha^{ -1} \beta \|}A^{\|\alpha^{ -1} 1_n\|} B^{\|\beta^{ -1} 1_n\|} \ .
\end{align}
\end{definition} 
Now we show that $M(X,Y,A,B)$ and $K(X,Y,A,B)$ are related by the boolean transform $\mathcal F_{\Boole}$.
\begin{theorem}\label{theorem:compatible transform} 
We have
\eq{
\mathcal F_{\Boole} :K(X,Y,A,B)\mapsto M(X,Y,A,B)
}
\end{theorem}
\begin{proof}
Following the notations in Definition \ref{definition:MK} and using Lemma \ref{lemma:decomposition},
\eq{
m_n
&=  \sum_{\sigma \in \inter(n) }\sum_{\substack{\alpha \in \inter(n) \\ \beta \in L (n) \\ \alpha \vee_{\inter} \beta = \sigma}} 
Y^{\|\alpha^{-1} \beta\|}A^{\|\alpha \|} B^{\|\beta\|} 
&=   \sum_{\sigma \in \inter(n) } \prod_{c \in \sigma } 
\sum_{\substack{\alpha \in \inter(|c|) \\ \beta \in L (|c|) \\ \alpha \vee_{\inter} \beta = 1_{|c|}}} 
Y^{\|\alpha^{-1} \beta\|}A^{\|\alpha \|} B^{\|\beta \|} \ .
}
Now, taking Kreweras complement, we have
\eq{
\sum_{\substack{\alpha \in \inter (|c|) \\ \beta \in L(|c|) \\ \alpha \vee_{\inter} \beta = 1_m}} 
Y^{\|\alpha^{-1} \beta\|}A^{\|\alpha \|} B^{\|\beta\|}
= \sum_{\substack{\alpha \in \kr \inter(|c|) \\ \beta \in \kr L (|c|) \\ \alpha\wedge_{\kr \inter} \beta = 0_m}} 
Y^{\|\alpha^{ -1} \beta \|}A^{\|\alpha^{-1} 1_{|c|} \|} B^{\|\beta^{-1} 1_{|c|} \|} = \kappa_{|c|}\ ,
}
where we used \eqref{eq:kreweras_rewritten}, \eqref{eq:kreweras-loops} and \eqref{eq:meet-join-kr}.
Applying Definition \ref{definition:moment-cumulant2} completes the proof.
\end{proof}

\subsection{Useful lemmas}
In this subsection we collect claims to be used in the following sections. Readers can come back later when they are needed. 

\begin{figure}[htbp!]
\includegraphics[width=0.8\textwidth]{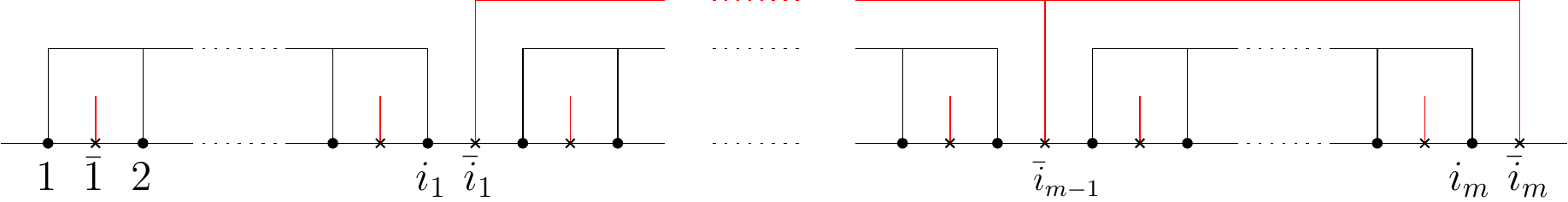}
\caption{Showing how to construct elements of $\kr \inter (n)$ where $n_m = n$. The support of the non-trivial red block is $Q \sqcup \{n\}$ in Lemma \ref{lemma:krint}.}
\label{figure:comb}
\end{figure}

\begin{lemma}\label{lemma:krint}
We have the following identification:
\eq{
\kr\inter(n)  = \{ Q \sqcup \{n\} : Q \subseteq  [n-1] \} \ .
}
\end{lemma}
\begin{proof}
Take some interval partition $\alpha = (1 \ldots, i_1) ( i_{1}+1 \ldots, i_2)   \ldots (i_{m-1} +1, \ldots, i_m) \in \inter(n)$ with $i_m= n$. Then, by the definition of Kreweras complement, the elements $i_1, i_2, \ldots,n$ constitute a block in the complement, but other elements are always isolated. See Figure \ref{figure:comb}.
\end{proof}

\begin{remark} \label{remark:krint}
We make some remarks on Lemma \ref{lemma:krint}.
\begin{enumerate}
\item The identification shows that each element in $\kr\inter(n)$ has at most one non-trivial block $Q \sqcup \{n\}$, which would necessarily contain $n$. We call this block a \emph{comb}.
\item The notation $Q_\alpha \subseteq [n-1]$ means that the subset $Q_\alpha$ is induced from $\alpha \in \kr\inter(n)$.
\end{enumerate}
\end{remark}

\begin{lemma}\label{lemma:wedge-int}
For $\alpha \in \kr\inter(n)$ and $\beta \in \nc(n)$, we have
\eq{
\alpha \wedge_{\kr \inter} \beta = (Q_{\alpha} \sqcup \{n\}) \cap \beta(n) \ ,
}
where $\beta(n)$ is a block of $\beta$ containing $n$.
In particular, if $\beta \in \kr\inter(n)$ then
\eq{
\alpha \wedge_{\kr \inter} \beta = \alpha \wedge \beta \ .
}
\end{lemma}
\begin{proof}
By using Lemma \ref{lemma:krint} and Definition \ref{definition:vee-int}, we have
\eq{
\alpha \wedge_{\kr\inter} \beta =\max \{ Q \sqcup 
\{n\}  : Q \subseteq  [n-1], \, Q \sqcup \{n\} \leq Q_\alpha \sqcup \{n\}, \beta(n) \}= (Q_\alpha \sqcup \{n\} ) \cap \beta(n) \ .
}
Next, $\beta \in \kr\inter(n)$ implies
\eq{
\alpha \wedge_{\kr \inter} \beta  = (Q_\alpha \cap Q_\beta) \sqcup \{n\} = \alpha \wedge \beta \ .
}
This completes the proof.
\end{proof}

\begin{figure}[htbp!]
\begin{subfigure}{.49\linewidth}
\centering
\includegraphics[width=0.75\textwidth]{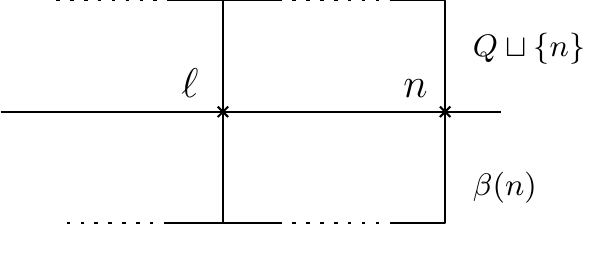}
\caption{$\ell \in Q_{\alpha}  \cap \beta(n) \not = \emptyset$ contradicts $\alpha \wedge_{\kr\inter} \beta = 0_n$.}
\label{figure:case0}
\end{subfigure}
\begin{subfigure}{.49\linewidth}
\centering
\includegraphics[width=0.9\textwidth]{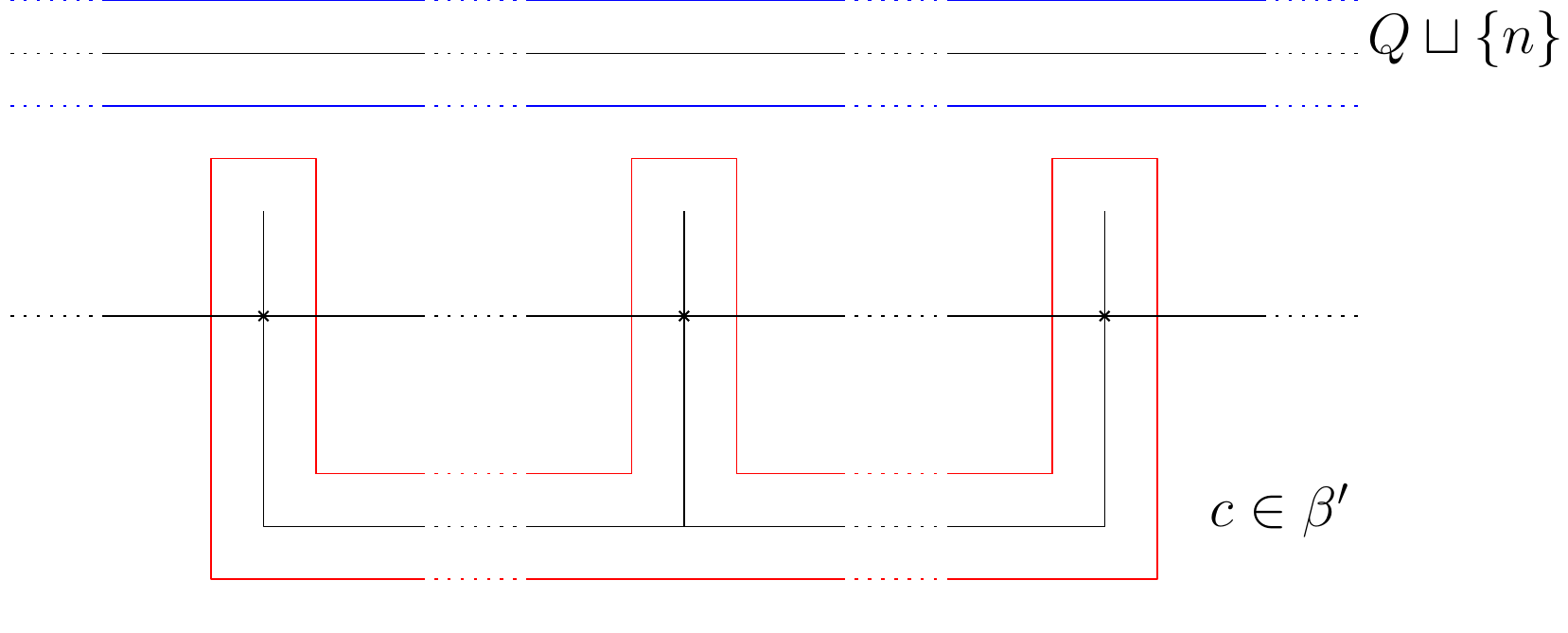}
\caption{One new loop is made if $|Q\cap c| = 0$.}
\label{figure:case1}
\end{subfigure}\\[3mm]
\begin{subfigure}{.49\linewidth}
\centering
\includegraphics[width=0.9\textwidth]{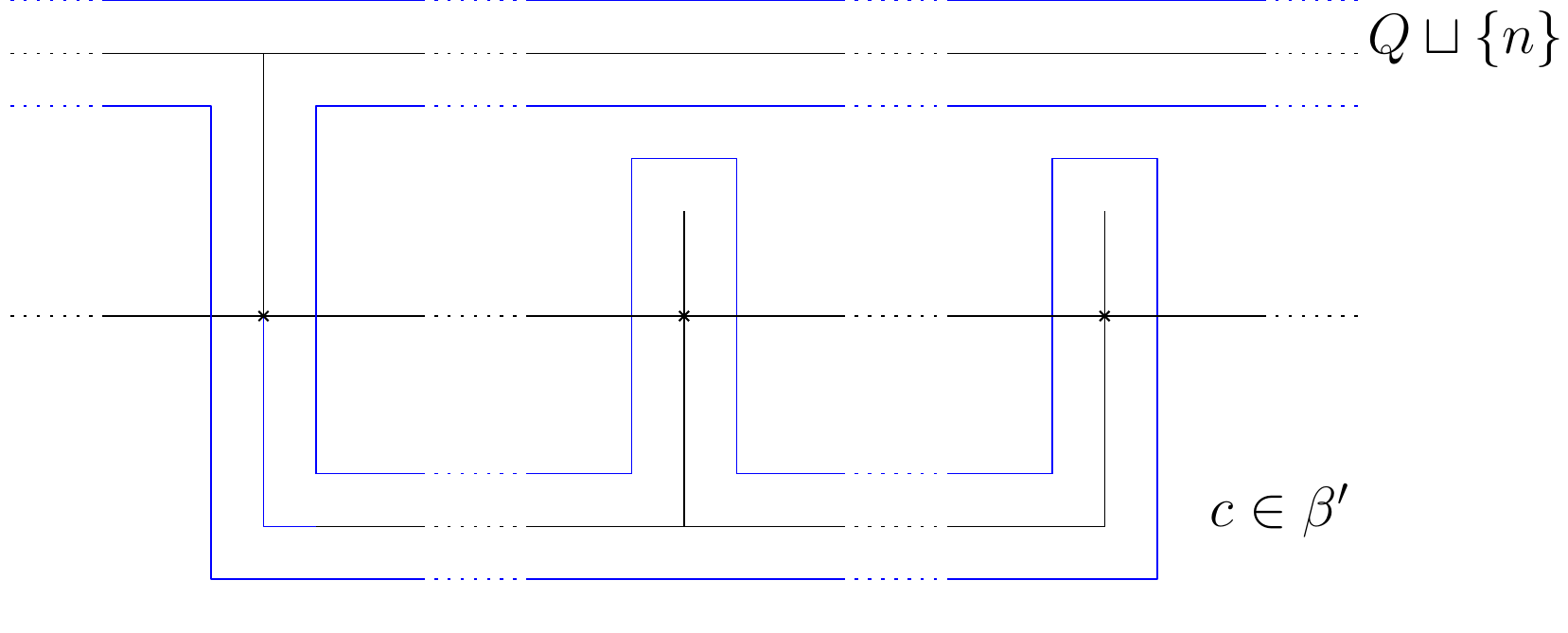}
\caption{No new loop will be made if $|Q\cap c| =1$.}
\label{figure:case2}
\end{subfigure}
\begin{subfigure}{.49\linewidth}
\centering
\includegraphics[width=0.9\textwidth]{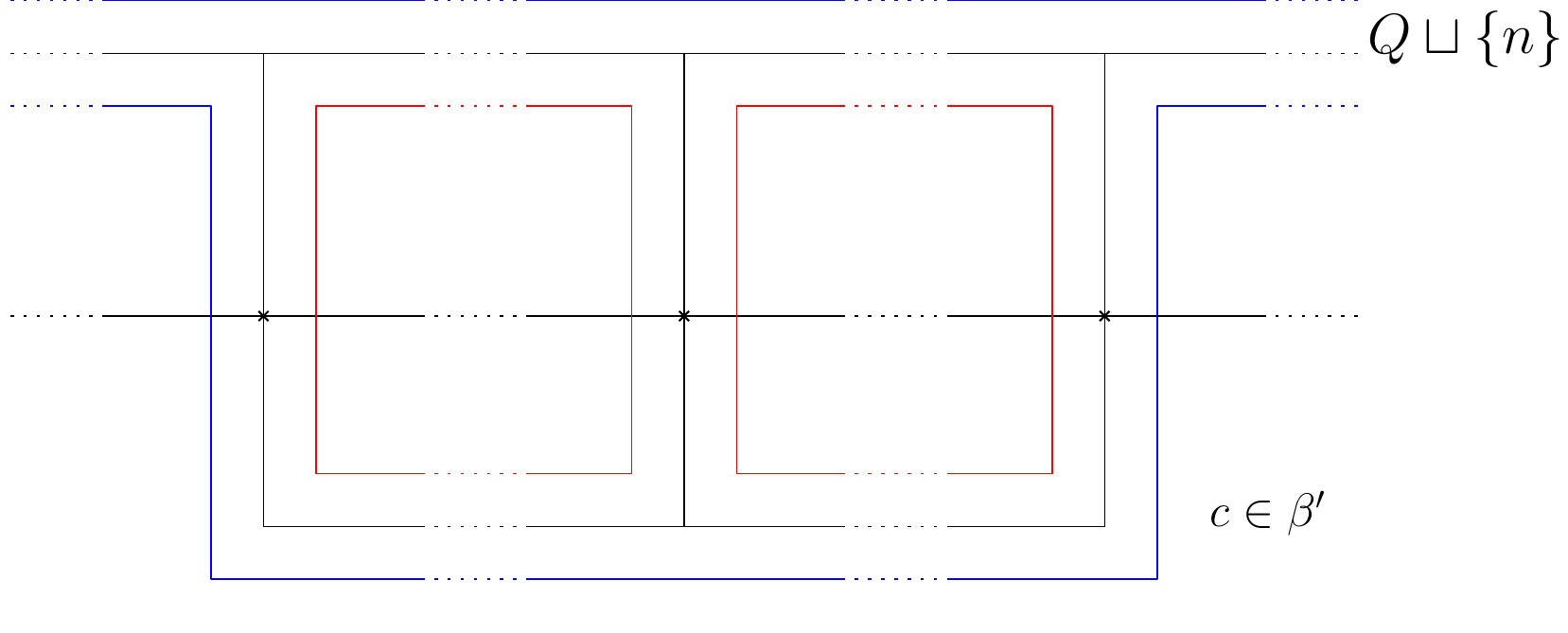}
\caption{$(|Q\cap c|-1)$ new loops ``inside'' if $|Q\cap c| \geq 2$.}
\label{figure:case3}
\end{subfigure}
\caption{Figures \ref{figure:case1}, \ref{figure:case2}, \ref{figure:case3} describe \eqref{eq:3cases} with $|c|=3$. The upper black lines represent $Q\sqcup \{n\}$ and the lower black lines $c$. Preexisting loops of meandric systems are represented by blue lines and new loops produced by adding $c$ are indicated by red lines.
}
\label{fig:test}
\end{figure}

\begin{lemma}\label{lemma:counting-loops}
For $\alpha \in \kr\inter (n)$ and $\beta \in \nc(n)$ with $\alpha \wedge_{\kr\inter} \beta = 0_n$, we have
\eq{
\#(\alpha^{ -1} \beta ) = 2 \cdot \# \{ c \in \beta^\prime: Q \cap c = \emptyset  \} +1 - \#(\beta^\prime) +  |Q|
}
where $\beta = \beta^\prime \oplus \beta(n)$ such that $\beta(n)$ is a block containing $n$, and $Q = Q_\alpha$.
\end{lemma}
\begin{proof}
We count loops in the meandric system made of $(\alpha,\beta)$ by adding cycles in $\beta$ one by one.
First, the condition $\alpha \wedge_{\kr\inter} \beta = 0_n$ and Lemma \ref{lemma:wedge-int} imply that $Q  \cap \beta(n) = \emptyset$.
Figure \ref{figure:case0} shows that having $\ell \in Q \cap \beta(n)$ 
would contradict the condition $\alpha \wedge_{\kr\inter} \beta = 0_n$.
This means that adding $\beta(n)$ does not increase the number of loops, which corresponds to case 2 below. 
Next, we add a cycle $c \in \beta^\prime$ to increase the number of loops as follows:
\eq{\label{eq:3cases}
\text{case 1: }\quad |Q\cap c| = 0  &\Rightarrow  +1 \\
\text{case 2: }\quad|Q\cap c| = 1 &\Rightarrow   \pm 0 \\
\text{case 3: }\quad|Q\cap c|  \geq 2  &\Rightarrow  + |Q\cap c|  -1 \\
}
Although case 3 includes case 2, 
we use the above classification to make things clear. 
First the condition $ |Q\cap c| = 0$ means that the block $c$ produces a loop without any interaction with $Q\sqcup \{n\}$ as in Figure \ref{figure:case1}, where a newly created loop is drawn in red. Second, with the condition $ |Q\cap c| = 1$, no new loop will be created although the preexisting loop containing $n$ will be stretched by $c$, which is drawn by the blue line in Figure \ref{figure:case2}. Third, in case $|Q\cap c|  = m \geq 2$, suppose $Q\cap c = \{i_1, \ldots, i_m \}$ and we connect $Q\sqcup \{n\}$ and $c$ one after another. To begin with, $i_1$ does not make any loop as in the case 2. Next, however, connection at $i_2$ gives a new loop, which will be enclosed by the preexisting loop,  increasing genus by one; see Figure \ref{figure:case3}. 
This inductive argument shows the claim on case 3.

Therefore, 
\eq{
\#(\alpha^{ -1} \beta )  = 1+  \sum_{c \in \beta^\prime}  \left[ 2 \cdot 1_{Q\cap c = \emptyset}  -1 + |Q\cap c| \right]
}
which leads to the formula because $Q \cap \beta(n) = \emptyset$ implies
\eq{
\sum_{c \in \beta^\prime} |Q \cap c| = |Q| \ .
}
This completes the proof.
\end{proof}

\begin{lemma}\label{lemma:binomial-formula}
For any partition $\beta$ of order $n$, not necessarily in $\nc(n)$, 
\eq{
\sum_{Q \subseteq [m]} A^{|Q|} B^{\# \{c\in\beta: Q \cap c = \emptyset\}} = \prod_{c \in \beta} \left( (A+1)^{|c|} +B-1\right)
}
\end{lemma}
\begin{proof}
First note that the LHS, which we denote by $S(\beta)$, is multiplicative for block decomposition in $\nc (n)$. Indeed, for $\beta = \beta_1 \oplus \beta_2$
\eq{
S(\beta) &= \sum_{Q \subseteq [n]} A^{|Q|} B^{\# \{c\in\beta_1: Q \cap c = \emptyset\}} B^{\# \{c\in\beta_2: Q \cap c = \emptyset\}} \\
 &= \sum_{Q_1 \subseteq [\beta_1]}  \sum_{Q_2 \subseteq [\beta_2]}A^{|Q_1|+|Q_2|} B^{\# \{c\in\beta_1: Q_1 \cap c = \emptyset\}} B^{\# \{c\in\beta_2: Q_2 \cap c = \emptyset\}} = S(\beta_1) S(\beta_2)
}
where $[\beta_i]$ is the support of $\beta_i$. 
Clearly the RHS is also multiplicative, so we prove the formula only for the case $\beta= 1_m$. Indeed,
\eq{
(\text{LHS}) = \sum_{Q \subseteq [m]} A^{|Q|} -1 +B = (A+1)^m  +B -1 = (\text{RHS})
}
where we treated the case $Q = \emptyset$ separately.
\end{proof}

\section{Thin meandric systems}\label{sec:thin}
In this section we consider the case where paths on the both upper and lower sides of the coordinate line consist of interval partitions, i.e.~$(\alpha, \beta) \in \inter(n) \times \inter(n)$. Since such a meandric system has both a shallow top and a shallow bottom (using the terminology from \cite{goulden2020asymptotics}), we shall call them \emph{thin meandric systems}. Since in this case there is no complicated layer structure due to non-crossing partitions, all the calculations are straightforward. 

\begin{theorem}
For meandric systems of $\inter(n) \times \inter(n)$,
the moment generating function $M(X,Y,A,B)$ and the cumulant generating function $K(X,Y,A,B) $ in Definition \ref{definition:MK}
are calculated as follows. 
\eq{
M(X,Y,A,B) &= \frac {X}{1-X(1+AB + (A+B)Y)} \quad \text{ and} \\
K(X,Y,A,B) &= \frac{X}{1-X(AB + (A+B)Y)} \ .
}
\end{theorem}
\begin{proof}
Since the generating function $M(X,Y,A,B)$ is obtained from $K(X,Y,A,B)$ by using Proposition \ref{proposition:Boole}:
\eq{
M(X,Y,A,B) &=  \frac{X}{1-X(AB+(A+B)Y)} \left[ 1- \frac{X}{1-X(AB+(A+B)Y)} \right]^{-1} \\
&= \frac {X}{1-X(1+AB + (A+B)Y)} \ .
}
we calculate $K(X,Y,A,B)$ for the rest of the proof.
\eq{
K(X,Y,A,B) 
&=\sum_{n=1}^\infty X^n 
\sum_{\substack{\alpha, \beta \in \kr\inter(n) \\ \alpha  \wedge_{\inter} \beta = 0_n}} 
Y^{\|\alpha^{ -1} \beta \|}A^{\|\alpha^{ -1} 1_{n} \|} B^{\|\beta^{ -1} 1_{n} \|} \\
&= \sum_{n=1}^\infty X^n 
\sum_{\substack{Q, R \subseteq [n-1]\\ Q \cap R = \emptyset}} 
Y^{|Q|+|R|}A^{n-1-|Q|} B^{n-1-|R|} \ .
}
Here, $Q\sqcup \{n\}$ and $R\sqcup \{n\}$ are the supports of the combs (see Remark \ref{remark:krint}) of $\alpha$ and $\beta$, and clearly
\eq{
\alpha \wedge_{\inter} \beta = 0_n \quad \Longleftrightarrow \quad Q \cap R = \emptyset \ ,
}
which also implies $\|\alpha^{ -1} \beta \| = |Q|+|R|$.
In addition, \eqref{eq:kreweras_rewritten} and \eqref{eq:kreweras_norm} explain the powers of $A$ and $B$.

Therefore,
\eq{
K(X,Y,A,B) &=
\sum_{n=1}^\infty X^n 
\sum_{\substack{Q, R \subseteq [n-1]\\ Q \cap R = \emptyset}} 
(AB)^{n-1-|Q|-|R|} (AY)^{|R|} (BY)^{|Q|} \\
&= \sum_{n=1}^\infty X^n 
\left( AB  + \left(A + B\right) Y \right)^{n-1} 
= \frac{X}{1-X(AB+(A+B)Y)} \ .
}
This completes the proof. 
\end{proof}

\begin{corollary}\label{cor:thin-number-loops}
The number of thin meandric systems of order $n$ having $k$ connected components is given by
$\displaystyle  2^{n-1}\binom{n-1}{k-1}$.
\end{corollary}
\begin{proof}
Set $A=B=1$ in the above result and extract the coefficient of  $X^nY^{n-k}$:
$$
[X^n Y^{n-k}] M(X,Y) = [X^n Y^{n-k}] \left(\frac{X}{1-2X(1+Y)}\right) = [X^n Y^{n-k}] X \sum_{n=0}^\infty (2X(1+Y))^n \ .
$$
This completes the proof.
\end{proof}
Notice that thin meanders (i.e.~$k=1$ above) correspond to $Q$ and $R$ forming a partition of $[n-1]$, hence there are $2^{n-1}$ such objects.

\section{Meandric systems with shallow top}\label{sec:shallow-top}

This section contains one of the main results of the paper, a generating series for the number of meandric systems with shallow top. The terminology comes from \cite{goulden2020asymptotics}, where meandric systems having one partition (say, the top one) being an interval partition have been called \emph{shallow top meanders}. It was recognized in \cite{goulden2020asymptotics} that this restricted setting allows for an explicit enumeration of meanders (meandric systems with one connected component), due to the simpler structure of the arches involved. 

We derive the cumulant generating function of shallow top meandric systems, and then apply the machinery from Section \ref{sec:transforms} to obtain the moment generating function. Our results generalize \cite[Theorem 1.1]{goulden2020asymptotics} adding two new statistics to the generating function: the number of loops of the meandric system (counted by $Y$) and the number of cycles of the non-interval partition ($\beta$ in our notation, counted by $B$).

\begin{theorem}\label{thm:ST}
For meandric systems of $\inter(n) \times \nc(n)$,
the boolean cumulant generating function $K(X,Y,A,B) $ in Definition \ref{definition:MK}
is given by
\begin{equation}\label{eq:K-ST}
K(X) = h(X(1+\hat g(X))).
\end{equation}
Here, $\hat g = \mathcal F_{\free} (g)$, and $g(x)$ and $h(x)$ are defined as
\eq{\label{eq:gh}
g(X) &= \sum_{n=1}^\infty g_n X^n &\text{where} \quad&g_n=  BY\left[ (1+AY)^n + (AY)^n (Y^{-2} -1)  \right] \ , \\
h(X) &= \sum_{n=1}^\infty h_n X^n & \text{where} \quad& h_n = (AY)^{n-1} \ .
}
The generating function for shallow top meandric systems is given by
$$M(X,Y,A,B) = \sum_{n=1}^\infty X^n \sum_{\substack{\alpha \in \inter(n) \\ \beta \in NC(n)}} 
Y^{\|\alpha^{-1} \beta\|}A^{\|\alpha\|} B^{\|\beta\|} =  \frac{K(X,Y,A,B)}{1-K(X,Y,A,B)}.$$
\end{theorem}
\begin{proof}
Using Lemma \ref{lemma:counting-loops} with the decomposition $\beta = \beta^\prime \oplus \beta(n)$,
\eq{
&K(X,Y,A,B) \\
&=\sum_{n=1}^\infty X^n \sum_{ \beta \in \nc (n) } 
\sum_{\substack{\alpha \in \kr\inter(n) \\ \alpha \wedge_{\inter} \beta = 0_n}} 
Y^{\|\alpha^{-1} \beta \|}A^{\|\alpha^{-1} 1_{n} \|} B^{\|\beta^{-1} 1_{n} \|} \\
&= \sum_{n=1}^\infty X^n \sum_{ \beta \in \nc (n) } 
\sum_{Q \subseteq [\beta^\prime]}  
Y^{ n-  2 \cdot \# \{ c \in \beta^\prime: Q \cap c = \emptyset  \} -1 + \#(\beta^\prime) - |Q|} 
A^{n-1-|Q|} B^{ \#(\beta^\prime)} \ .
}
Here, $[\beta^\prime]$ is the support of $\beta^\prime$ and $Q = Q_{\alpha}$,
and $\|\beta^{-1} 1_{n} \| =n-1- \|\beta\|  = n-1- (n- \#(\beta)) = \#(\beta^\prime)$.
Moreover, we used the fact that
\eq{
\alpha \wedge_{\inter} \beta = 0_n \Longleftrightarrow Q \cap \beta(n) = \emptyset \ .
}
Then, we continue our calculation with Lemma \ref{lemma:binomial-formula}:
\eq{
&K(X,Y,A,B) \\&= \sum_{n=1}^\infty X^n \sum_{ \beta \in \nc (n) } 
(AY)^{n-1}
(BY)^{ \#(\beta^\prime) } 
\sum_{Q \subseteq [\beta^\prime]}  \left( A^{-1}Y^{-1} \right)^{|Q|}
\left(Y^{ -  2} \right)^{\# \{ c \in \beta^\prime: Q \cap c = \emptyset  \} } \\
&= \sum_{n=1}^\infty X^n \sum_{ \beta \in \nc (n) } 
(AY)^{|\beta^\prime|+|\beta (n)|-1}
(BY)^{ \#(\beta^\prime) } 
\prod_{c \in \beta^\prime} \left[ (A^{-1}Y^{-1} +1)^{|c|}+Y^{-2}-1   \right] \\
&= \sum_{n=1}^\infty X^n \sum_{ \beta \in \nc (n) } 
(AY)^{|\beta(n)|-1}
\prod_{c \in \beta^\prime} BY\left[ (1+AY)^{|c|} + (AY)^{|c|} (Y^{-2} -1)  \right]
}
Therefore, using the definition \eqref{eq:gh} and Lemma \ref{lemma:free2} we calculate
\eq{
K(X,Y,A,B)
= \sum_{n=1}^\infty X^n \sum_{ \beta \in \nc (n) } 
h_{|\beta(n)|} \prod_{c \in \beta^\prime} g_{|c|} 
= \sum_{s=1}^\infty    h_{s} X^s
\left( 1+ \sum_{i=1}^\infty \hat g_i X^i \right)^s \ .
}
where $\mathcal F_\free: g \mapsto \hat g$ is the free transform. This completes the proof. 
\end{proof}

Note that the boolean cumulant generating function $K$ from Eq.~\eqref{eq:K-ST} is not quite explicit, due to the definition of the function $\hat g$, which is given implicitly through its free transform $\mathcal F_{\free}$. Solving for $\hat g$ explicitly requires inverting a function, which cannot be done in full generality. Our theorem has the theoretical interest of expressing the moment generating function of the shallow top meandric systems with the help of the functional transforms associated to to the type of lattices the top, respectively the bottom partitions belong to. Moreover, one can use the implicit formulas from Theorem \ref{thm:ST} to extract useful information regarding the enumeration of shallow top meandric systems. For example, the main result of \cite{goulden2020asymptotics} states that the number of shallow top meanders on $2n$ points with $m$ blocks on the bottom is 
$$\frac 1 n \binom{n}{m-1}\binom{n+m-1}{m-1} = [X^n Y^{n-1} A^{n-m}]M(X,Y,A,1).$$

\section{Shallow top semi-meandric systems}\label{sec:shallow-top-semi}

In this section we consider \emph{shallow top semi-meandric systems}: meandric systems formed by
interval partitions and so-called \emph{rainbow partition}, which is defined as follows:
\eq{
\beta_{\mathrm rainbow} =  \begin{cases} \displaystyle (1, n) (2,n-1) \cdots (\frac n2, \frac n2 +1) & \text{($n$ is even)}\\[10pt] \displaystyle (1, n) (2,n-1) \cdots (\frac{n+1}2) & \text{($n$ is odd)} \ . \end{cases}
}
The terminology is justified by the bijection between the set of semi-meandric systems and the set of meandric systems where one of the partitions (say the bottom one) is fixed to be the rainbow partition \cite{difrancesco1997meander}. We further specialize this setting by considering interval partitions on the top. 

\begin{figure}[htbp!]
\begin{subfigure}{.49\linewidth}
\centering
\includegraphics[width=0.7\textwidth]{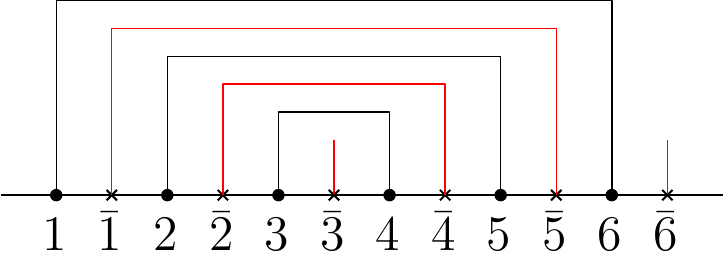}
\caption{For $\beta_{\mathrm rainbow} = (1,6)(2,5)(3,4)$, we have\\
$\beta_{\mathrm rainbow}^{\kr} = (\bar 1, \bar 5)(\bar 2, \bar 4)(\bar 3)(\bar 6)$.}
\label{figure:rainbow_even}
\end{subfigure}
\begin{subfigure}{.49\linewidth}
\centering
\includegraphics[width=0.7\textwidth]{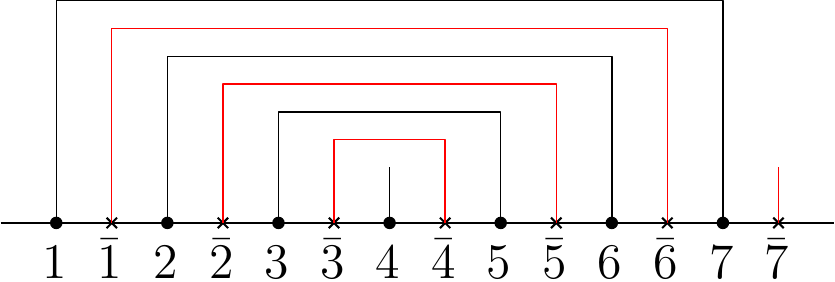}
\caption{For $\beta_{\mathrm {rainbow}} = (1,7)(2,6)(3,5)(4)$, we have\\
$\beta_{\mathrm rainbow}^{\kr} = (\bar 1, \bar 6)(\bar 2, \bar 5)(\bar 3, \bar4)(\bar 7)$.}
\label{figure:rainbow_odd}
\end{subfigure}
\caption{Examples for the rainbow partition $\beta_{\mathrm rainbow}$ when $n=6,7$, and its Kreweras complement.}
\label{figure:rainbow}
\end{figure}

One can find graphical representations of rainbow partitions (as well as their Kreweras complements) in Figure \ref{figure:rainbow}. By using the decomposition in \eqref{eq:decomposition}, we can write
\begin{equation}\label{eq:rainbow-Kr-n}
    \beta_{\mathrm rainbow}^{\kr} = (\beta_{\mathrm rainbow}^{\kr})^\prime
+ \beta_{\mathrm rainbow}^{\kr}(n).
\end{equation}
It always holds that $\beta_{\mathrm rainbow}^{\kr}(n) = \{n\}$, and abusing the notation, $(\beta_{\mathrm rainbow}^{\kr})^\prime \in NC(n-1)$. Moreover, all cycles in $(\beta_{\mathrm rainbow}^{\kr})^\prime$ are of length 2, unless $n$ is even, when only one exceptional cycle consists of a single point: $\{n/2\}$.

Note that the definition of the moment generating function in Definition \ref{definition:MK} can be naturally extended to the current case, so we can state and prove the main result of this section. 

\begin{theorem}\label{theorem:rainbow}
The generating function of shallow top semi-meandric systems 
$M(X,Y,A) $ defined by $\inter(n) \times \{\beta_{\mathrm rainbow}\}$ is given by
\begin{equation}\label{eq:m-rainbow}
    M(X,Y,A)  = \frac{ X + X^2(Y+A)}{ 1 - X^2 Y(1+2YA + A^2)}.    
\end{equation}
 
\end{theorem}

\begin{proof}
In this proof, we treat our problem in the Kreweras-complement view, but do not use the cumulant generating function. First, by writing $\beta := \beta_{\mathrm rainbow}^{\kr}$, we have
\eq{
M(X,Y,A) &= \sum_{n=1}^\infty X^n \ 
\sum_{\alpha \in \inter(n)} Y^{\|\alpha^{-1}\beta_{\mathrm rainbow} \|} A^{\|\alpha\|}\\
& = \sum_{n=1}^\infty X^n \ 
\sum_{\substack{\alpha \in \kr\inter(n) \\ \alpha \wedge_{\inter} \beta = 0_n}} Y^{\|\alpha^{-1}\beta \|} A^{\|\alpha^{-1}1_n\|} \ , 
}
where we used Proposition \ref{proposition:kreweras_loop}.
We claim that the condition $\alpha \wedge_{\inter} \beta = 0_n$ always holds in this case. Indeed, as in Remark \ref{remark:krint}, $\alpha \in \kr\inter(n)$ consists of isolated points and possibly at most one non-trivial cycle containing $n$. On the other hand, we see from 
Eq.~ \eqref{eq:rainbow-Kr-n} 
and Figure \ref{figure:rainbow}
that $\{n\}$ is always an isolated point in $\beta$. Hence, Definition \ref{definition:vee-int} implies that the meet of the two partitions inside the lattice $\kr \inter(n)$ is trivial.  

Next, we apply Lemma \ref{lemma:counting-loops} to $NC(n)$ and then Lemma \ref{lemma:binomial-formula} to $NC(n-1)$.
\eq{
M(X,Y,A) 
& = \sum_{n=1}^\infty X^n \sum_{Q \subseteq [n-1]} \
Y^{n-  2 \cdot \# \{ c \in \beta^\prime: Q \cap c = \emptyset  \} -1 +\lfloor \frac n 2 \rfloor - |Q|}
A^{n-1-|Q|}  \\
& =  \sum_{n=1}^\infty X^n 
Y^{n-1 +\lfloor \frac n 2 \rfloor }A^{n-1}
\sum_{Q \subseteq [n-1]} 
(Y^{-2})^{ \# \{ c \in \beta^\prime: Q \cap c = \emptyset  \} }
((YA)^{-1})^{|Q|}  \\
& =  \sum_{n=1}^\infty X^n 
\
(YA)^{n-1} Y^{\lfloor \frac n 2 \rfloor } 
\underbrace{\prod_{c \in \beta^\prime}
\left(  ((YA)^{-1} +1)^{|c|}  + Y^{-2} -1  \right)}_{(\star)} \ .
}
Moreover, recall that the cycles of $\beta'$ are all pairs, except for the case when $n$ is even, where we have an extra singleton. Hence, we can make $(\star)$ explicit:
\eq{
(\star) &= \left(  (YA)^{-2} +2(YA)^{-1}   + Y^{-2}  \right)^{\left\lfloor \frac{n}{2} \right\rfloor -e(n)}
\left( (YA)^{-1}   + Y^{-2}  \right)^{e(n)} 
}
where 
\eq{
e(n) = \begin{cases}
1 & \text{ $n$ is even}\\
0 & \text{ $n$ is odd}.
\end{cases}
}
Here, note that $\displaystyle 2\left(\left\lfloor \frac{n}{2} \right\rfloor -e(n)\right) + e(n) = n-1$, which is the number of points in $\beta'$.

Finally, dividing the series into two parts depending on the parity of $n$, we have
\eq{
M(X,Y,A) 
& =  \sum_{n=1}^\infty X^n 
\
[Y(1+2YA + A^2)]^{\left\lfloor \frac{n}{2} \right\rfloor -e(n)} (Y+A)^{e(n)}\\
& = X \sum_{m=0}^\infty [X^2 Y(1+2YA + A^2)]^m
+ X^2(Y+A)\sum_{m=0}^\infty [X^2 Y(1+2YA + A^2)]^m \\
& = \frac{ X + X^2(Y+A)}{ 1 - X^2 Y(1+2YA + A^2)} \ .
}
This completes the proof.
\end{proof}

\begin{remark}
It is straightforward to extract the distribution of the number of loops at fixed $n$ from Eq.~\eqref{eq:m-rainbow}:
$$\forall n \geq 1, \qquad [X^n]M(X,Y,1) = \begin{cases}
(2Y)^{k-1}Y^k &\qquad\text{if }n=2k\\
(2Y(Y+1))^{k-1} &\qquad\text{if }n=2k-1.
\end{cases}$$
In particular, the number of shallow top semi-meanders is $2^{\lceil n/2\rceil-1}$.
\end{remark}

\section{Random matrix models for meandric systems}\label{sec:RMT}

We present in this section several matrix models for the different types of meandric systems that we study. These models are motivated by quantum information theory and they allow for a uniform presentation, the matrix model being constructed from a tensor product of two (random) completely positive maps related respectively to the type of non-crossing partitions used to build the meander.

\subsection{Meanders}

We start with the case of usual meanders and meandric systems, obtained by stacking two general non-crossing partitions one on top of the other. We shall first state two models from the literature and then introduce a new, simpler one, which will be generalized in later subsections to different types of meandric systems. 

Let us define the \emph{meander polynomial}
$$m_n(\ell) = \sum_{\alpha, \beta \in NC(n)} \ell^{\#(\alpha^{-1}\beta)},$$
where $\alpha,\beta$ are the non-crossing partitions used to build the meandric system having $\#(\alpha^{-1}\beta)$ loops (or connected components). Note that $m_n$ is the coefficient of $X^n$ of the following polynomial $M$ (similar to the one in  \eqref{eq:def-M}, see also \cite[Section 5]{fukuda2019enumerating}), evaluated at $A=B=1$
\begin{equation}\label{eq:def-M-NC-NC}
M(X,Y,A,B) = \sum_{n=1}^\infty X^n \sum_{\alpha,\beta \in NC (n)} 
Y^{\|\alpha^{-1} \beta\|}A^{\|\alpha \|} B^{\|\beta\|} 
\end{equation}

The first matrix models for meandric systems are due to P.~Di Francesco and collaborators, see \cite[Section 5]{difrancesco1997meander} or \cite[Section 6]{di2001matrix}. We recall here, for the sake of comparison with our new models, the GUE-based construction from the former reference. A \emph{Ginibre random matrix} $G$ is a matrix having independent and identically distributed (i.i.d.) entries $G_{ij}$ following a standard complex Gaussian distribution; a Ginibre random matrix can be rectangular, and we do not assume any symmetry properties for it. Define now a \emph{GUE} (Gaussian Unitary Ensemble) random matrix 
$$B = \frac{G+G^*}{\sqrt 2} \in \mathcal M_d(\mathbb C)^{sa},$$
where $G$ is a $d \times d$ Ginibre matrix. Note that the GUE matrix defined above is not normalized in the usual way, see \cite[Chapter 2]{anderson2010introduction} or \cite[Chapter 1]{mingo2017free}. 

\begin{proposition}\cite[Section 5]{difrancesco1997meander}
Let $\ell$ be a fixed positive integer and consider $B_1, \ldots, B_\ell \in \mathcal M_d(\mathbb C)$ i.i.d.~GUE matrices. Then, for all $n \geq 1$, 
$$m_n(\ell) = \lim_{d \to \infty} \mathbb E \frac{1}{d^2} \operatorname{Tr} \left( \sum_{i=1}^\ell\frac{B_i \otimes \bar B_i}{d} \right)^{2n}.$$
\end{proposition}

A second matrix model for meanders was discovered in relation to the theory of quantum information, more precisely in the study of partial transposition of random quantum states. We recall briefly the setup here. A \emph{Wishart random matrix} of parameters $(d,s)$ is simply defined by $W=GG^*$, where 
$G \in \mathcal M_{d \times s}(\mathbb C)$ 
is a Ginibre matrix. Note that $W$ is by definition a positive semidefinite matrix, and thus its normalized version $\rho = W / \operatorname{Tr} W$ is called a \emph{density matrix} in quantum theory \cite{watrous2018theory}. This model for random density matrices was introduced in \cite{zyczkowski2001induced} and it is called the \emph{induced measure} of parameters $(d,s)$. For bi-partite quantum states $\rho \in \mathcal M_{d^2}(\mathbb C) = \mathcal M_{d}(\mathbb C) \otimes \mathcal M_{d}(\mathbb C)$, the \emph{partial transposition} operation 
$$ \rho^\Gamma := [\operatorname{id}_d \otimes \operatorname{transp}_d] (\rho)$$
plays a crucial role in quantum information theory, in relation to the notion of entanglement \cite{horodecki2009quantum}. Before stating the result from \cite{fukuda2013partial}, let us also mention that the combinatorics of meanders appears also in computations related to random quantum channels, see \cite[Section 6.2]{fukuda2015additivity}. 
\begin{proposition}\cite[Theorem 4.2]{fukuda2013partial}
Let $\rho \in \mathcal M_{d^2}(\mathbb C)$ be a random bi-partite quantum state of parameters $(d^2, \ell)$ for some fixed integer $\ell \geq 1$. Then, for all $n \geq 1$, 
$$m_n(\ell) = \lim_{d \to \infty} \mathbb E \frac{1}{d^2} \operatorname{Tr} \left( \ell d \rho^\Gamma \right)^{2n}.$$
\end{proposition}

We would like to introduce now a new, simpler, matrix model for meandric systems, which we shall later generalize to include different types of non-crossing partitions. To begin, recall the classical Stinespring dilation theorem \cite{stinespring1955positive} from operator theory: any completely positive (CP) map $\Phi : \mathcal M_{d}(\mathbb C) \to \mathcal M_{d'}(\mathbb C)$ can be written as 
$$\Phi(X) =  [\operatorname{id}_{d'} \otimes \operatorname{Tr}_s](AXA^*),$$
for some operator $A : \mathbb C^d \to \mathbb C^{d'} \otimes \mathbb C^s$; taking $s = dd'$ allows one to recover all CP maps by varying $A$. We shall denote the map above by $\Phi_A$.  Moreover, by imposing the condition that $A$ is an isometry, one obtains in this way all completely positive and trace preserving maps, i.e. all \emph{quantum channels} \cite[Corollary 2.27]{watrous2018theory}. We also introduce the (un-normalized) maximally entangled state 
$$\Omega_d = \sum_{i=1}^d e_i \otimes e_i \in \mathbb C^d \otimes \mathbb C^d$$
for a fixed basis $\{e_i\}_{i=1}^d$ of $\mathbb C^d$ and the rank-one matrix (having trace $d$)
$$\omega_d = \Omega_d \Omega_d^* \in \mathcal M_{d^2}(\mathbb C).$$

\begin{theorem}\label{thm:RM-model-NC-NC}
Consider two independent Ginibre matrices $G,H \in \mathcal M_{\ell \times  d^2}(\mathbb C)$ and the corresponding CP maps $\Phi_{G,H}:\mathcal M_\ell(\mathbb C) \to \mathcal M_d(\mathbb C)$. Define
$$Z := [\Phi_G \otimes \Phi_H](\omega_\ell) \in \mathcal M_{d^2}(\mathbb C)$$
Then, for all $n \geq 1$, 
$$m_n(\ell) = \lim_{d \to \infty} \mathbb E \frac{1}{d^2} \operatorname{Tr} \left(\frac{Z}{d^2}\right)^n.$$
\end{theorem}
\begin{proof}
The statement is a moment computation which is quite standard in the theory of random matrices. We give a proof using the graphical version of Wick's formula developed in \cite{collins2011gaussianization}. The diagram corresponding to the matrix $Z$ is depicted in Figure \ref{fig:Z}. 

\begin{figure}[htbp!]
\centering
\includegraphics[scale=1]{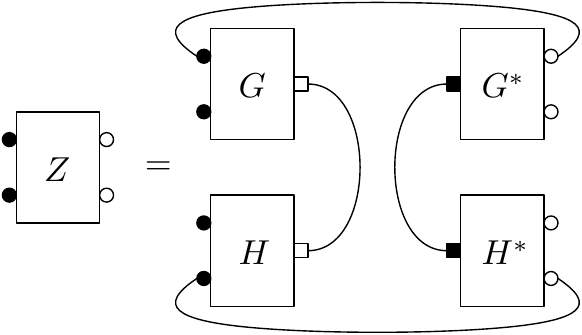}
\caption{A graphical representation of the bi-partite matrix $Z = [\Phi_G \otimes \Phi_H](\omega_\ell)$. Round labels represent the space $\mathbb C^d$, while square labels represent $\mathbb C^\ell$.}
\label{fig:Z}
\end{figure}

To compute the $n$-th moment of $Z$, $\mathbb E \operatorname{Tr} Z^n$, one considers the expectation value of the trace of the concatenation of $n$ instances of the diagram in Figure \ref{fig:Z}. This expectation value is, according to the graphical Wick formula \cite[Theorem 3.2]{collins2011gaussianization}, a combinatorial sum indexed by two permutations $\alpha, \beta \in \mathcal S_n$ of diagrams $\mathcal D_{\alpha, \beta}$. Note that we have here two (independent) permutations since we are dealing with independent Gaussian matrices $G$ and $H$: the permutation $\alpha$ is encoding the wiring of the $G$-matrices, while $\beta$ encodes the wiring of the $H$-matrices. A diagram $\mathcal D_{\alpha, \beta}$ consists of (see Figure \ref{fig:E-Tr-Z-2} for a simple example): 
\begin{itemize}
\item $\# \alpha + \#(\alpha^{-1}\gamma)$ loops corresponding to the round decorations of $G$
\item $\# \beta + \#(\beta^{-1}\gamma)$ loops corresponding to the round decorations of $H$
\item $\#(\alpha^{-1}\beta)$ loops corresponding to all the square decorations,
\end{itemize}
\begin{figure}[htbp!]
\centering
\includegraphics[scale=1]{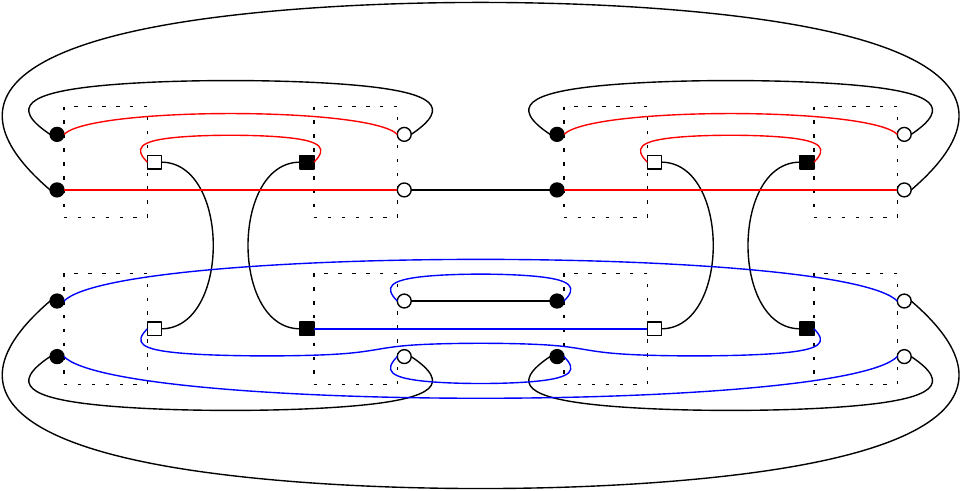}
\caption{The diagram $\mathcal D_{\textcolor{red}{\alpha=(1)(2)},  \textcolor{blue}{\beta=(12)}}$ as a term in the graphical Wick expansion of $\mathbb E \operatorname{Tr} Z^2$. This diagram consists of 6 $d$-dimensional loops (3 corresponding to the top (red) pictures and 3 corresponding to the bottom (blue) ones) and of 1 $\ell$-dimensional loop.}
\label{fig:E-Tr-Z-2}
\end{figure}
where $\gamma = (1 \, 2 \, 3 \, \cdots \, p) \in \mathcal S_p$ is the full cycle permutation and the boxes are numbered from right to left. Let us first justify the formula for the number of $d$-dimensional (i.e.~corresponding to round decorations) loops given by $\alpha$. Note that the permutation $\alpha$, acting on the top boxes, gives rise to two types of loops: the top ones, in which the output of the $i$-th $G$ box is connected to the corresponding input of the $i$-th $G^*$ box, and the bottom ones, where the output of the $i$-th $G$ box is connected to the corresponding input of the $\gamma(i)$-th $G^*$ box. 
It is now an easy combinatorial fact that the number of connected components of a bipartite graph on $2n$ vertices having edges 
$$\{(i, n+\gamma(i))\}_{i=1}^n \sqcup \{(i, n+\sigma(i))\}_{i=1}^n$$
is precisely $\#(\sigma^{-1}\gamma)$, proving our claim. A similar argument settles the case of the $\ell$-dimensional loops, where the top and the bottom symbols are identified by the wires corresponding to the maximally entangled state $\omega_d$. 

The result of applying the graphical Wick formula is thus
$$d^{-2-2p}\mathbb E \operatorname{Tr} Z^p = d^{-2-2p} \sum_{\alpha,\beta \in \mathcal S_n} d^{\# \alpha + \#(\alpha^{-1}\gamma) + \# \beta + \#(\beta^{-1}\gamma)} \ell^{\#(\alpha^{-1}\beta)}.$$
Standard combinatorial inequalities about permutations (\cite{biane1997some} or \cite[Lecture 23]{nica2006lectures}) give
\begin{align*}
\#\alpha + \#(\alpha^{-1}\gamma) &\leq p+1 \\
\#\beta + \#(\beta^{-1}\gamma) &\leq p+1
\end{align*}
with equality if and only if both $\alpha$ and $\beta$ are geodesic permutations (see Section \ref{sec:NC-permutations}) corresponding to non-crossing partitions. Moreover, for such permutations, $\#(\alpha^{-1}\beta)$ is precisely the number of loops of the meandric system built from $\alpha$ and $\beta$ (see Proposition \ref{proposition:loops}), finishing the proof. 
\end{proof}

\begin{remark}
In the statement above, one can replace the random CP map $\Phi_H$ with $\Phi_G$ or even $\Phi_{\bar G}$. This fact, quite surprising at first, is due to the particular asymptotic regime we are interested in, that is $d \to \infty$ and $\ell$ fixed. When performing the Gaussian integration using the graphical Wick calculus, one obtains a sum over permutations $\alpha \in \mathcal S_{2n}$; however, due to the fact that $\ell$ is fixed, the permutation $\alpha$ will be constraint to leave invariant the top (resp.~the bottom) $n$ points; this, in turn, amounts to having a decomposition $\alpha = \alpha^T \oplus \alpha^B$, with $\alpha^{T,B} \in \mathcal S_n$, and the proof would continue as above. Note that if $\ell$ would grow with $d$, different behavior would occur, see e.g.~\cite{collins2010random}.
\end{remark}

\begin{remark}
One can keep track of the parameters $A$ and $B$ appearing in the definition of the generating function $M$ from eq.~\eqref{eq:def-M-NC-NC} by adding a decoration of type ``A'' (resp.~``B'') on the partial traces appearing in the Stinespring dilation formulas for the channels $\Phi_G$ (resp.~$\Phi_H$); we leave the details to the reader. 
\end{remark}

\subsection{Shallow top meanders}

We consider in this section shallow top meanders, that is meanders built out of a general non-crossing partition and an interval partition (which sits on the top). We shall construct a random matrix model for these combinatorial objects, by replacing the random channel $\Phi_H$ from Theorem \ref{thm:RM-model-NC-NC} by a non-random channel. First, we define the corresponding shallow top meander polynomial by
$$m_n^{\mathrm{ST}}(\ell):=\sum_{\substack{\alpha \in \inter(n)\\ \beta \in \nc(n)}}  \ell^{\#(\alpha^{-1}\beta)}=\sum_{\substack{\alpha \in \kr\inter(n)\\ \beta \in \nc(n)}} \ell^{\#(\alpha^{-1}\beta)}.$$

\begin{theorem}\label{thm:RM-model-Int-NC}
Consider a Ginibre matrix $G \in \mathcal M_{\ell \times  d^2}(\mathbb C)$ and the corresponding completely positive map 
\begin{align*}
    \Phi_{G}:\mathcal M_\ell(\mathbb C) &\to \mathcal M_d(\mathbb C)\\
    X &\mapsto [\trace_d \otimes \id_d](GXG^*).
\end{align*}
Define
\begin{align*}
Z &:= [\Phi_G \otimes \Psi](\omega_\ell) \in \mathcal M_{d\ell}(\mathbb C)\\
Z_0 &:= [\Phi_G \otimes \id](\omega_\ell) \in \mathcal M_{d\ell}(\mathbb C),
\end{align*}
where the completely positive map $\Psi$ is defined by
\begin{align}
    \label{eq:def-Psi} \Psi:\mathcal M_\ell(\mathbb C) &\to \mathcal M_\ell(\mathbb C)\\
    \nonumber X &\mapsto X + (\trace X) I_\ell.
\end{align}
Then, for all integers $n,\ell \geq 1$, 
$$m_n^\mathrm{ST}(\ell) = \lim_{d \to \infty} \mathbb E \frac{1}{d} \operatorname{Tr}[(d^{-1}Z_0) (d^{-1}Z)^{n-1}] .$$
\end{theorem}
\begin{proof}
We shall use the graphical Wick formula to compute the expectation value $\mathbb E \trace[Z_0Z^{n-1}]$. We shall encode the action of the linear map $\Psi$ by its Choi-Jamio{\l}kowski matrix $C_\Psi = \omega_{\ell} + I_{\ell^2}$. Diagrammatically, we shall apply the Wick formula to the diagram in Figure \ref{fig:ST-Wick}, with
$$\quad C_1 = C_2 = \cdots = C_{n-1} = C_\Psi = \omega_\ell + I_{\ell^2} \quad \text{ and } C_n = \omega_\ell.$$
\begin{figure}[htbp!]
    \centering
    \includegraphics{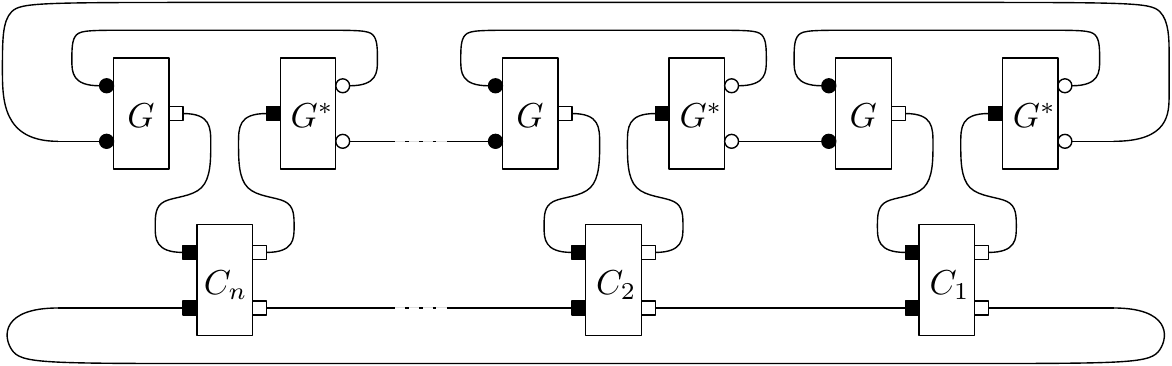}
    \caption{The diagram corresponding to 
    $\trace[Z_0Z^{n-1}]$.}
    \label{fig:ST-Wick}
\end{figure}

Applying the graphical Wick formula to compute the expectation over the Gaussian random matrix $G$, we have 
$$\mathbb E \trace[Z_0Z^{n-1}] = \sum_{\beta \in \mathcal S_n} d^{\#\beta}d^{\#(\beta^{-1}\gamma)} \trace_{\beta,\gamma}[C_1, C_2, \ldots, C_n],$$
where the trace factor above corresponds to the diagram obtained by connecting the top output of the $i$-th $C$-box to the top input of the $\beta(i)$-th $C$-box, and the bottom output of the $i$-th $C$-box to the bottom input of the $\gamma(i)$-th $C$-box, see Figure \ref{fig:ST-C-trace}. 
\begin{figure}[htbp!]
    \centering
    \includegraphics{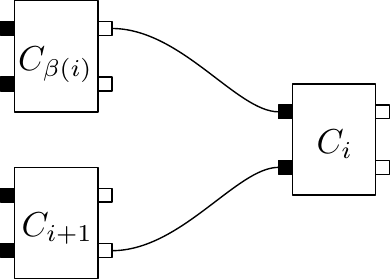}
    \caption{Connecting the Choi-Jamio{\l}kowski matrices $C_i$ by the permutation $\beta$ and $\gamma$, where $\gamma(i) = i+1$.}
    \label{fig:ST-C-trace}
\end{figure}
Note that the matrices $C_i$ are of finite size $\ell^2$. Thus, in order to take the limit $d \to \infty$, we have to maximize the exponent $\#\beta + \#(\beta^{-1}\gamma)$. Using the triangle inequality, we obtain (see the proof of Theorem \ref{thm:RM-model-NC-NC})
$$\lim_{d \to \infty} \mathbb E \frac{1}{d} \operatorname{Tr}[(d^{-1}Z_0) (d^{-1}Z)^{n-1}] =
\sum_{\beta \in \nc(n)} \trace_{\beta,\gamma}[C_1, C_2, \ldots, C_n].$$
We shall now develop the diagram corresponding to the trace in the sum above. We shall encode the choice of $\omega_\ell$ or $I_{\ell^2}$ for each matrix $C_i$ (here, $i \in [n-1]$) by a subset $Q \subseteq [n-1]$: an integer $i \in [n-1]$ is an element of $Q$ if and only if we choose the matrix $\omega_\ell$ for the box $C_i$. Let $\alpha \in \kr \inter(n)$ be the comb partition encoded by the subset $Q$ (see Lemma \ref{lemma:krint}). We claim that
\begin{equation}\label{eq:claim-ST-loops}
    \trace_{\beta,\gamma}[C^Q_1, C^Q_2, \ldots, C^Q_{n-1}, C_n] = \ell^{\#(\alpha^{-1}\beta)},
\end{equation}
where, for $i \in [n-1]$,
$$C^Q_i = \begin{cases}
\omega_\ell &\qquad \text{ if } i \in Q\\
I_{\ell^2} &\qquad \text{ if } i \notin Q.
\end{cases}$$
The claim \eqref{eq:claim-ST-loops} allows us to conclude, since
$$m^{\mathrm{ST}}_n(\ell) = \sum_{\substack{Q \subseteq [n-1] \\ \beta \in \nc(n)}} \trace_{\beta,\gamma}[C^Q_1, C^Q_2, \ldots, C^Q_{n-1}, C_n] = \sum_{\substack{\alpha \in \kr \inter(n)\\ \beta \in \nc(n)}}\ell^{\#(\alpha^{-1}\beta)}.$$

Let us now prove \eqref{eq:claim-ST-loops}. Recall from Lemma \ref{lemma:krint} that the geodesic comb permutation $\alpha \in \kr \inter(n)$ associated to a subset $Q \subseteq [n-1]$ is given by
$$\alpha(i) = \begin{cases}
\bar q_{j+1} &\qquad \text{ if } i = \bar q_j \in \bar Q\\
i &\qquad \text{ if } i\notin \bar Q,
\end{cases}$$
where $\bar Q = Q \sqcup \{n\} = \{\bar q_1, \ldots, \bar q_{|Q|+1} \}$. Hence, if we were to follow the top outputs of the boxes $C^Q_i$, we would have (see Figure \ref{fig:ST-Q}):
$$ i \mapsto \begin{cases}
\beta_i = \alpha^{-1}\circ \beta(i) &\qquad \text{ if } \beta(i) \notin \bar Q\\
\bar q_{j-1} = \alpha^{-1}\circ \beta(i) &\qquad \text{ if } \beta(i) = \bar q_j \in \bar Q,
\end{cases}$$
which shows the claim \eqref{eq:claim-ST-loops}, finishing the proof. 
\begin{figure}[htbp!]
    \centering
    \includegraphics{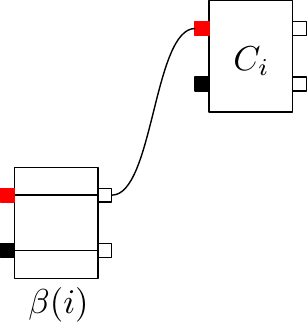} \qquad\qquad\qquad\qquad \includegraphics{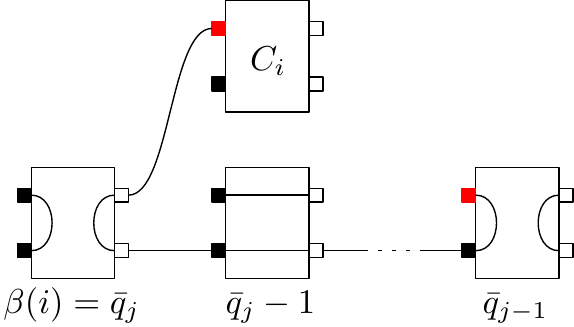}
    \caption{Following the top outputs (in red) of the $C$-boxes in the diagram $\trace_{\beta,\gamma}[C^Q_1, C^Q_2, \ldots, C^Q_{n-1}, C_n]$. Left: $\beta(i) \notin \bar Q$; right: $\beta(i) = \bar q_j \in \bar Q$.}
    \label{fig:ST-Q}
\end{figure}
\end{proof}

\subsection{Thin meandric systems}

In the case of thin meandric systems (corresponding to bottom and top permutations corresponding to interval partitions, see Section \ref{sec:thin}), there is a matrix model which is closely related to the one in the previous section. Actually, one needs to replace in the statement of Theorem \ref{thm:RM-model-Int-NC} the random CP map $\Phi_G$ (responsible for the general non-crossing permutation $\beta$) by another copy of the deterministic linear CP map $\Psi$ from \eqref{eq:def-Psi}. 
Before stating and proving the result, let us define the corresponding meander polynomial:
$$m_n^{\mathrm{thin}}(\ell):=\sum_{\alpha,\beta \in \inter(n)}  \ell^{\#(\alpha^{-1}\beta)}=\sum_{\alpha,\beta \in \kr\inter(n)} \ell^{\#(\alpha^{-1}\beta)}.$$

\begin{theorem}\label{thm:RM-model-Int-Int}
Recall the linear, completely positive map $\Psi$ from \eqref{eq:def-Psi} and define the matrix $Z := [\Psi \otimes \Psi](\omega_\ell) \in \mathcal M_{d^2}(\mathbb C)$. Then, for all integers $n,\ell \geq 1$, 
$$m_n^\mathrm{thin}(\ell) =  \operatorname{Tr}[\omega_l Z^{n-1}] = \ell(2+2\ell)^{n-1}.$$
\end{theorem}
\begin{proof}
First, Figure \ref{fig:cj} shows how we can interpret the Choi-Jamio{\l}kowski matrix of $\Psi$.
\begin{figure}[htbp!]
    \centering
    \includegraphics{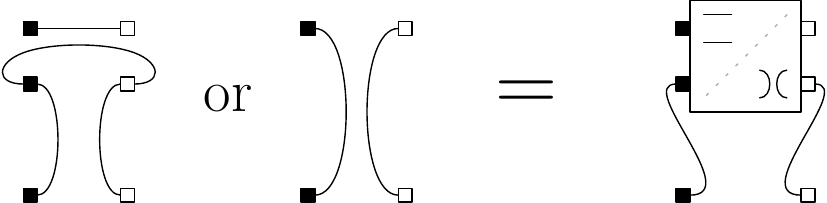}
    \caption{The left hand side show the Choi-Jamio{\l}kowski matrices of the two possible operations; $X \mapsto \trace[X]I_\ell$ or $X \mapsto X$, respectively.}
    \label{fig:cj}
\end{figure}

Then, it is straightforward to see that the diagram corresponding to $m_n^\mathrm{thin}(\ell)$ is the one from Figure \ref{fig:thin}, 
where there are $(n-1)$ boxes containing the sum of $I_{\ell^2}$ and
$\omega_\ell$ on each of the two rows. 
\begin{figure}[htbp!]
    \centering
    \includegraphics{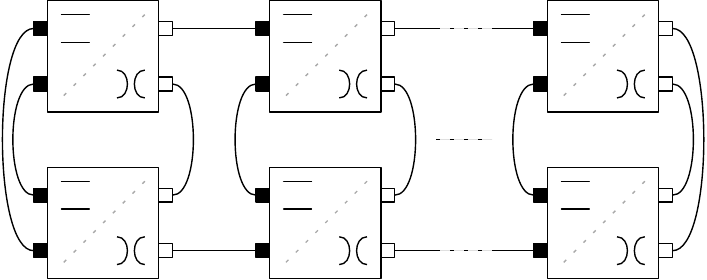}
    \caption{The diagram for $\operatorname{Tr}[\omega_l Z^{n-1}]$. Each row contains $(n-1)$ boxes corresponding to the Choi-Jamio{\l}kowski matrix of the map $\Psi$.}
    \label{fig:thin}
\end{figure}
Develop now the diagram as a sum indexed by pairs $(Q,R)$ of subsets of $[n-1]$, where in the top row we replace the $i$-th box by $\omega_\ell$ if $i \in Q$ and by the identity matrix otherwise, and we use the subset $R$ in the similar manner for the bottom row. It is straightforward to see that the diagram obtained has at most $n$ loops (each contributing a factor $\ell$), and that the exact number of loops is 
$n - |Q \Delta R|$, where $\Delta$ is the symmetric difference operation. In other words, each time the $i$-th boxes are different on the two rows, a loop is ``lost''. Hence, 
$$\operatorname{Tr}[\omega_l Z^{n-1}] = \sum_{Q,R \subseteq [n-1]} \ell^{n-|Q\Delta R|}.$$
It is now easy to check that, given two permutations $\alpha, \beta \in \kr \inter(n)$ defined respectively by the subsets $Q,R \subseteq [n-1]$, we have 
$$\#(\alpha^{-1}\beta) = n - |Q \Delta R|,$$
establishing the first claim. The final equality is obtained by noting but
$$Z = (2+\ell)I + \omega_l = (2+\ell)\left(I - \frac{\omega_l}{l}\right) + (2+2l) \frac{\omega_l}{l},$$
which implies in turn that
$$Z^{n-1} =  (2+\ell)^{n-1}\left(I - \frac{\omega_l}{l}\right) + (2+2l)^{n-1} \frac{\omega_l}{l}$$
and thus $\operatorname{Tr}[\omega_l Z^{n-1}] = \ell (2+2l)^{n-1}$; see also Corollary \ref{cor:thin-number-loops}.
\end{proof}

\section*{Acknowledgment}
MF acknowledges JSPS KAKENHI Grant Number JP16K00005. 
This work was supported by Bilateral Joint Research Projects (JSPS, Grant number JPJSBP120203202 and MEAE-MESRI, PHC Sakura).

\bibliography{ref}{}
\bibliographystyle{alpha}

\end{document}